\documentclass[11pt]{amsart}
\usepackage{amsmath, amssymb, latexsym}
\usepackage{mathrsfs}
\usepackage[all, knot]{xy}
\xyoption{arc}

\theoremstyle{definition}

\numberwithin{equation}{section}

\numberwithin{equation}{section}

\newtheorem{theorem}{Theorem}[section]
\newtheorem{corollary}[theorem]{Corollary}
\newtheorem{lemma}[theorem]{Lemma}
\newtheorem{proposition}[theorem]{Proposition}

\theoremstyle{definition}
\newtheorem{definition}[theorem]{Definition}
\newtheorem{remark}[theorem]{Remark}

\newcommand{\HOM}{\text{HOM}}
\newcommand{\Hom}{\text{Hom}}

\newcommand{\B}{\mathbb{B}}

\newcommand{\Q}{\mathbb{Q}}
\newcommand{\R}{\mathscr{R}}
\newcommand{\Z}{\mathbb{Z}}
\newcommand{\N}{\mathbb{N}}

\newcommand{\g}{\mathfrak{g}}
\newcommand{\K}{\mathbb{K}}

\newcommand{\ii}{\textit{\textbf{i}}}
\newcommand{\jj}{\textit{\textbf{j}}}
\newcommand{\kk}{\textit{\textbf{k}}}
\newcommand{\Seq}{\text{Seq}}

\newcommand{\Pol}{\mathscr{P}}

\newcommand{\Ind}{\text{Ind}}
\newcommand{\Res}{\text{Res}}
\newcommand{\Span}{\text{Span}}
\newcommand{\gdim}{\mathbf{Dim}}
\newcommand{\Mod}{\text{Mod}}
\newcommand{\pMod}{\text{pMod}}
\newcommand{\fMod}{\text{fMod}}
\newcommand{\pr}{\text{pr}}
\newcommand{\infl}{\text{infl}}

\newcommand{\hd}{\text{hd }}
\newcommand{\soc}{\text{soc }}
\newcommand{\Max}{\text{max}}
\newcommand{\re}{\text{re}}
\newcommand{\im}{\text{im}}
\newcommand{\iso}{\text{iso}}
\newcommand{\Indd}{\widehat{\text{Ind}}}

\newcommand{\genO}[3] % generator 1
{
\fontsize{9}{9}\selectfont
\xy
(0,5)*{}; (0,-5)*{} **\dir{-};
(4,0)*{\cdots};
(8,5)*{}; (8,-5)*{} **\dir{-};
(12,0)*{\cdots};
(16,5)*{}; (16,-5)*{} **\dir{-};
(0,-7)*{#1}; (8,-7)*{#2}; (16,-7)*{#3};
\endxy
\fontsize{10}{10}\selectfont}
\newcommand{\genX}[3] % generator x
{
\fontsize{9}{9}\selectfont
\xy
(0,5)*{}; (0,-5)*{} **\dir{-};
(4,0)*{\cdots};
(8,5)*{}; (8,-5)*{} **\dir{-};
(8,0)*{\bullet}; (12,0)*{\cdots};
(16,5)*{}; (16,-5)*{} **\dir{-};
(0,-7)*{#1}; (8,-7)*{#2}; (16,-7)*{#3};
\endxy
\fontsize{10}{10}\selectfont}
\newcommand{\genT}[4] % generator tau
{\fontsize{9}{9}\selectfont
\xy
(0,5)*{}; (0,-5)*{} **\dir{-};
(4,0)*{\cdots};
(5,5)*{}; (12,-5)*{} **\dir{-};
(12,5)*{}; (5,-5)*{} **\dir{-};
(13,0)*{\cdots};
(17,5)*{}; (17,-5)*{} **\dir{-};
(0,-7)*{#1}; (5,-7)*{#2}; (12,-7)*{#3}; (17,-7)*{#4};
\endxy
\fontsize{10}{10}\selectfont}

\newcommand{\dcross}[2]
{\fontsize{9}{9}\selectfont
\xy
(0,7)*{}="T1"; (7,7)*{}="T2";
(0,-7)*{}="B1"; (7,-7)*{}="B2";
"T1"; "B1" **\crv{(11, 0)};
"T2"; "B2" **\crv{(-4,0)};
(0,-9)*{#1}; (7,-9)*{#2};
\endxy
\fontsize{10}{10}\selectfont}
\newcommand{\dcrossA}[2]
{\fontsize{9}{9}\selectfont
\xy
(0,5)*{}; (0,-5)*{} **\dir{-};
(6,5)*{}; (6,-5)*{} **\dir{-};
(0,-7)*{#1}; (6,-7)*{#2};
\endxy
\fontsize{10}{10}\selectfont}
\newcommand{\dcrossL}[3]
{ \fontsize{10}{10}\selectfont
\xy
(0,5)*{}; (0,-5)*{} **\dir{-};
(6,5)*{}; (6,-5)*{} **\dir{-};
(0,0)*{\bullet}; (-5,0)*{#1}; (0,-7)*{#2}; (6,-7)*{#3};
\endxy
\fontsize{10}{10}\selectfont}
\newcommand{\dcrossR}[3]
{\fontsize{10}{10}\selectfont
\xy
(0,5)*{}; (0,-5)*{} **\dir{-};
(6,5)*{}; (6,-5)*{} **\dir{-};
(6,0)*{\bullet}; (11,0)*{#1}; (0,-7)*{#2}; (6,-7)*{#3};
\endxy
\fontsize{10}{10}\selectfont}

\newcommand{\LU}[3]
{\fontsize{9}{9}\selectfont
\xy
(0,5)*{}; (8,-5)*{} **\dir{-} \POS?(.25)="x";
(8,5)*{}; (0,-5)*{} **\dir{-};
"x"*{\bullet}; "x"+(-2,2)*{#1}; (0,-7)*{#2}; (8,-7)*{#3}; \endxy
\fontsize{10}{10}\selectfont}
\newcommand{\RD}[3]
{\fontsize{9}{9}\selectfont
\xy
(0,5)*{}; (8,-5)*{} **\dir{-} \POS?(.75)="x";
(8,5)*{}; (0,-5)*{} **\dir{-};
"x"*{\bullet}; "x"+(2,2)*{#1}; (0,-7)*{#2}; (8,-7)*{#3}; \endxy
\fontsize{10}{10}\selectfont}
\newcommand{\RU}[3]
{\fontsize{9}{9}\selectfont
\xy
(0,5)*{}; (8,-5)*{} **\dir{-};
(8,5)*{}; (0,-5)*{} **\dir{-} \POS?(.25)="x";
"x"*{\bullet}; "x"+(2,2)*{#1}; (0,-7)*{#2}; (8,-7)*{#3}; \endxy
\fontsize{10}{10}\selectfont}
\newcommand{\LD}[3]
{\fontsize{9}{9}\selectfont
\xy
(0,5)*{}; (8,-5)*{} **\dir{-};
(8,5)*{}; (0,-5)*{} **\dir{-} \POS?(.75)="x";
"x"*{\bullet}; "x"+(-2,2)*{#1}; (0,-7)*{#2}; (8,-7)*{#3}; \endxy
\fontsize{10}{10}\selectfont}

\newcommand{\BraidL}[3]
{\fontsize{9}{9}\selectfont
\xy
(0,7)*{}="T1"; (6,7)*{}="T2";  (12,7)*{}="T3";
(0,-7)*{}="B1"; (6,-7)*{}="B2"; (12,-7)*{}="B3";
"T1"; "B3" **\crv{(0,0)&(12,0)}; "T3"; "B1" **\crv{(12,0)&(0,0)};
"T2"; "B2" **\crv{(6,6)&(3,4.5)&(0,1.5)&(0,-1.5)&(3,-4.5)&(6,-6)};
(0,-9)*{#1}; (6,-9)*{#2};(12,-9)*{#3};
\endxy
\fontsize{10}{10}\selectfont}
\newcommand{\BraidR}[3]
{\fontsize{9}{9}\selectfont
\xy
(0,7)*{}="T1"; (6,7)*{}="T2";  (12,7)*{}="T3";
(0,-7)*{}="B1"; (6,-7)*{}="B2"; (12,-7)*{}="B3";
"T1"; "B3" **\crv{(0,0)&(12,0)}; "T3"; "B1" **\crv{(12,0)&(0,0)};
"T2"; "B2" **\crv{(6,6)&(8,4.5)&(12,1.5)&(12,-1.5)&(8,-4.5)&(6,-6)};
(0,-9)*{#1}; (6,-9)*{#2};(12,-9)*{#3};
\endxy
\fontsize{10}{10}\selectfont}
\newcommand{\threeDotStrands}[3]
{ \fontsize{9}{9}\selectfont
\xy
(0,7)*{}; (0,-7)*{} **\dir{-};
(5,7)*{}; (5,-7)*{} **\dir{-};
(10,7)*{}; (10,-7)*{} **\dir{-};
(0,0)*{\bullet};(10,0)*{\bullet}; (-3,0)*{c};(21,0)*{-a_{ij}-1-c}; (0,-9)*{#1}; (5,-9)*{#2};(10,-9)*{#3};
\endxy
\fontsize{10}{10}\selectfont}

\newcommand{\dcrossI}[2]
{\fontsize{9}{9}\selectfont
\xy
(0,6)*{}="T1"; (6,6)*{}="T2";
(0,-6)*{}="B1"; (6,-6)*{}="B2";
"T1"; "B1" **\crv{(10, 0)};
"T2"; "B2" **\crv{(-4,0)};
(0,-8)*{#1}; (6,-8)*{#2};
\endxy
\fontsize{10}{10}\selectfont}
\newcommand{\BraidLI}[3]
{\fontsize{9}{9}\selectfont
\xy
(0,6)*{}="T1"; (5,6)*{}="T2";  (10,6)*{}="T3";
(0,-6)*{}="B1"; (5,-6)*{}="B2"; (10,-6)*{}="B3";
"T1"; "B3" **\crv{(0,0)&(10,0)}; "T3"; "B1" **\crv{(10,0)&(0,0)};
"T2"; "B2" **\crv{(5,5)&(3,4)&(0,1.5)&(0,-1.5)&(3,-4)&(5,-5)};
(0,-8)*{#1}; (5,-8)*{#2};(10,-8)*{#3};
\endxy
\fontsize{10}{10}\selectfont}
\newcommand{\BraidRI}[3]
{\fontsize{9}{9}\selectfont
\xy
(0,6)*{}="T1"; (5,6)*{}="T2";  (10,6)*{}="T3";
(0,-6)*{}="B1"; (5,-6)*{}="B2"; (10,-6)*{}="B3";
"T1"; "B3" **\crv{(0,0)&(10,0)}; "T3"; "B1" **\crv{(10,0)&(0,0)};
"T2"; "B2" **\crv{(5,5)&(7,4)&(10,1.5)&(10,-1.5)&(7,-4)&(5,-5)};
(0,-8)*{#1}; (5,-8)*{#2};(10,-8)*{#3};
\endxy
\fontsize{10}{10}\selectfont}

\title[Quiver Hecke algebras for Borcherds-Cartan datum II]
{Quiver Hecke algebras for Borcherds-Cartan datum II}

\author[Bolun Tong]{Bolun Tong}
\address{Hankuk University of Foreign Studies, Seoul, Korea}
\email{tbl\_2018@hufs.ac.kr}

\author[Wan Wu]{Wan Wu${}^{*}$}
\address{Harbin Engineering University,
Harbin, China}
\email{wuwan1818@163.com}
\thanks{${}^{*}$ Corresponding author.}



\begin{document}
\maketitle
\begin{abstract}
We give the crystal structure of the Grothendieck group $G_0(R)$ of irreducible modules over the quiver Hecke algebra $R$ constructed in \cite{TW2023}. This leads to the categorification of the crystal $B(\infty)$ of the quantum Borcherds algebra $U_q(\g)$ and its irreducible highest weight crystal $B(\lambda)$ for arbitrary Borcherds-Cartan data. Additionally, we study the cyclotomic categorification of irreducible highest weight $U_q(\g)$-modules.
\end{abstract}

\section*{\textbf{Introduction}}

\vskip 2mm
In \cite{TW2023}, the authors introduced a family of quiver Hecke algebras (Khovanov-Lauda-Rouquier algebras) $R(\nu)$ for $\nu\in\N[I]$, associated with a Borcherds-Cartan datum $(I,A)$ consisting of a countable set $I$  and a symmetrizable Borcherds-Cartan matrix $A=(a_{ij})_{i,j\in I}$ indexed by $I$. Let $R=\bigoplus_{\nu\in \N[I]}R(\nu)$. Following the framework given by  Khovanov-Lauda \cite{KL2009,KL2011}, it was shown that the  Grothendieck group $K_0(R)$ of the category of finitely generated graded projective $R$-modules is isomorphic to the negative part of the corresponding quantum Borcherds algebra $U_q(\g)$. In particular, when $i$ is an imaginary index, i.e. $a_{ii}\leq 0$, the  algebra $R(ni)$ for $\nu=ni$ can be considered as a deformation of the nilcoxeter algebra.

 %$$\dcrossI{}{}=0\quad\quad \BraidLI{}{}{}  =  \BraidRI{}{}{}\quad\quad \LU{{}}{}{}    =  \RD{{}}{}{}\quad\quad \LD{{}}{}{}     =    \RU{{}}{}{} $$

In comparison with \cite{KOP2012}, which also provides a categorification of these algebras, our construction imposes no such restriction as  $a_{ii}\neq 0$. With this advantage, we could consider an arbitrary Borcherds-Cartan datum and use our construction to give a categorification of typical crystals and the irreducible highest weight modules over $U_q(\g)$.

Historically, in the Kac-Moody cases, A. Lauda and M. Vazirani \cite{LV2009} studied the crystal structure on the categories of graded modules over the quiver Hecke algebra $\mathcal R$, which categorify the negative part of the quantum Kac-Moody algebra $\mathcal U_q(\g)$ associated with a symmetrizable Cartan datum. They proved that, as crystals, the Grothendieck group of the graded irreducible modules over $\mathcal R$ is isomorphic to the crystal $\mathcal B(\infty)$ of the quantum Kac-Moody algebra, derived from Kashiwara's crystal basis theory. Let $\lambda\in P^+$ be a dominant integral weight. The crystal structure of the irreducible modules over the cyclotomic quiver Hecke algebra $\mathcal R^{\lambda}$ is isomorphic to the crystal $\mathcal B(\lambda)$ of the irreducible highest weight module $\mathcal V(\lambda)$. Consequently, one could compute the rank of the corresponding Grothendieck groups. In \cite{KK2012}, Kang and Kashiwara proved that the cyclotomic quotient $\mathcal R^{\lambda}$ provides a categorification of $\mathcal V(\lambda)$. That is, the Grothendieck group $\mathcal K_0(R^\lambda)$ of the finitely generated graded projective modules over $\mathcal R^\lambda$ has a natural $\mathcal U_q(g)$-module structure that is isomorphic to $\mathcal V(\lambda)$.

In the context of quantum Borcherds algebra, similar results have been obtained in \cite{KOP2012} and \cite{KKO2011}, with the condition that $a_{ii}\neq 0$. Particularly, in \cite{KOP2012}, a perfect basis theory for the quantum Borcherds algebra is developed to categorify the crystals.

We extend these results to an arbitrary Borcherds Cartan datum. We first define the Kashiwara operators on the Grothendieck group $G_0(R)$ of the irreducible modules, which is slightly different from the Kac-Moody cases. Then we categorify the crystal $B(\infty)$ of the quantum Borcherds algebra $U_q(\g)$ and its irreducible highest weight crystal $B(\lambda)$ for each $\lambda\in P^+$. In the final chapter, we construct an algebra $\mathscr R$ that is `larger' than $R$ and assert that the irreducible highest weight $U_q(\g)$-module $V(\lambda)$ can be realized as a subspace of $K_0(\R^\lambda)_{\Q(q)}$ or $G_0(\R^\lambda)_{\Q(q)}$ for the cyclotomic algebra $\R^\lambda$.

\vskip 3mm

\noindent\textbf{Notations.}
In this paper, $\K$ is a fixed algebraically closed field. For a graded $A$-module $M=\bigoplus_{n\in \Z}M_n$ over a $\Z$-graded $\K$-algebra  $A$, its graded dimension is defined to be
$$\gdim M=\sum_{n\in\Z}(\text{dim}_{\K}M_n)q^n,$$
where $q$ is a formal variable. For $m\in \Z$,  the degree shift $M\{m\}$ is the graded $A$-module obtained by setting $(M\{m\})_n=M_{n-m}$.  More generally, for  $f(q)=\sum_{m\in \Z}a_mq^m\in\N[q,q^{-1}]$, we set $M^{f}=\bigoplus_{m\in\Z}(M\{m\})^{\oplus a_m}$.
 \vskip 2mm
 Given two graded $A$-modules $M$ and $N$, $\text{Hom}_{A\text{-gr}}(M,N)$ is the ${\K}$-vector space of grading-preserving homomorphisms. We define the $\Z$-graded vector space $\text{HOM}_A(M,N)$ to be
$$\text{HOM}_A(M,N)=\bigoplus_{n\in\Z}\text{Hom}_{A\text{-gr}}(M\{n\},N)=\bigoplus_{n\in\Z}\text{Hom}_{A\text{-gr}}(M,N\{-n\}).$$

%We sometimes use the term gr-projective (resp. gr-irreducible, gr-free and so on) module in the sense of graded projective (resp. graded irreducible, graded free) module.

\vskip 6mm
\section{\textbf{Categorification of quantum Borcherds algebras}}

\subsection{Quantum Borcherds algebras}\

\vskip 2mm
Let $I$ be an index set possibly countably infinite. A symmetrizable Borcherds-Cartan matrix is an integer-valued matrix $A=(a_{ij})_{i,j\in I}$ satisfying
\begin{itemize}
\item[(i)] $a_{ii}=2, 0, -2, -4, \ldots$,

\item[(ii)] $a_{ij} \in \Z_{\le 0}$ for $i \neq j$,

\item[(iii)] there is a diagonal matrix $D=\text{diag} (r_i \in
\Z_{>0} \mid i \in I)$ such that $DA$ is symmetric.
\end{itemize}

\vskip 2mm
 We set $I^{\text{re}}=\{ i \in I \mid a_{ii}=2 \}$, $I^{\text{im}}= \{ i \in I \mid a_{ii} \le 0 \}$, and
$I^{\text{iso}} =\{ i \in I \mid a_{ii}=0 \}$.

\vskip 2mm
A Borcherds-Cartan datum is given by
 $$\begin{aligned}
 &  \text{a symmetrizable Borcherds-Cartan matrix}\ A,\\
 & \text{a  free abelian group}\ P,\ \text{called the weight lattice},\\
 & P^{\vee}=\Hom_\Z(P,\Z), \ \text{called the dual weight lattice},\\
& \alpha_i \in P  \ \text{for }  i \in I, \ \text{called the simple roots},\\
& h_i \in P^\vee  \ \text{for }  i \in I, \ \text{called the simple coroots},\\
& \text{a  symmetric bilinear form}\ ( \ , \ ):P\times P\rightarrow \Z.
\end{aligned} $$
such that $\alpha_j(h_i)=a_{ij}$ for all $i,j  \in I$, and $(\lambda,\alpha_i)=r_i\lambda(h_i)$ for all $\lambda \in P, i\in I$. In particular,
$$(\alpha_i,\alpha_j)=r_ia_{ij}=r_ja_{ji}.$$

The elements in $P^+=\{\lambda \in P \mid \lambda(h_i)\geq 0 \  \text{for all} \ i \in I \}$  are called  dominant integral weights. The free abelian group $Q=\bigoplus_{i \in I} {\Z \alpha_i}$ is called the root lattice. We identify the positive root lattice $Q_+=\sum_{i \in I}{\N\alpha_i}$ with $\N[I]$. %and write $\nu\cdot \nu'$ for $(\nu,\nu')$ $(\nu,\nu'\in \N[I])$ , so that
%$$i\cdot j=(\alpha_i,\alpha_j).$$

\vskip 2mm
Let $q$ be an indeterminate. For each $i\in I$, set $q_i=q^{r_i}$. For $i\in I^\re$ and $n\in \N$, we define
 $$[n]_i=\frac{q_i^n-q_i^{-n}}{q_i-q_i^{-1}}\ \text{and} \ [n]_i!=[n]_i[n-1]_i\cdots [1]_i.$$

The quantum Borcherds algebra $U(=U_q(\g))$ associated with a given Borcherds-Cartan datum is the $\Q(q)$-algebra  generated by $e_i,f_i\ (i\in I)$ and $q^h\ (h\in P^{\vee})$,  satisfying
$$\begin{aligned}
& q^0=1,\quad q^hq^{h'}=q^{h+h'} \quad\text{for} \ h,h' \in P^{\vee}, \\
& q^h e_{j}q^{-h} = q^{\alpha_j(h)} e_{j}, \ \ q^h f_{j}q^{-h} = q^{-\alpha_j(h)} f_{j}\quad \ h \in P^{\vee}, i\in I,\\
& e_if_j-f_je_i=\delta_{ij}\frac{K_i-K_i^{-1}}{q_i- q_i^{-1}} \quad\text{where} \ K_i=q^{r_ih_i},\\
& \sum_{r+s=1-a_{ij}}(-1)^r
{e_i}^{(r)}e_{j}e_i^{(s)}=0 \quad\text{for} \ i\in
I^{\text{re}},j\in I \ \text {and} \ i \neq j, \\
& \sum_{r+s=1-a_{ij}}(-1)^r
{f_i}^{(r)}f_{j}f_i^{(s)}=0 \quad\text{for} \ i\in
I^{\text{re}},j\in I \ \text {and} \ i \neq j, \\
& e_{i}e_{j}-e_{j}e_{i} = f_{i}f_{j}-f_{j}f_{i} =0 \quad\text{for}\ i,j\in I\ \text{and}\ a_{ij}=0.
\end{aligned}
$$
Here, $e_i^{(n)}=e_i^n/[n]_i!$ and $f_i^{(n)}=f_i^n/[n]_i!$ for $i\in I^\re$ and $n\in \N$. The algebra $U$ is $Q$-graded by assigning $|e_i|=\alpha_i$ and $|f_i|=-\alpha_i$.

\vskip 2mm

Let $U^{0}$ (resp. $U^+$, resp. $U^-$)  be the subalgebra of $U$ generated by $q^h$ for $h\in P^{\vee}$ (resp. $e_{i}$ for $i\in I$, resp. $f_{i}$ for $i\in I$).  Then we have the triangular decomposition
$U\cong U^-\otimes U^0 \otimes U^+$.
So for each $u\in U^-$, there exist unique $Q,R\in U^-$ such that
$$e_iu-ue_i=\frac{K_iQ-K_i^{-1}R}{q_i-q_i^{-1}}.$$
We set $e_i'(u)=R$. The operators $e_i'$ $(i\in I)$ satisfy the quantum Serre relations (cf. \cite[Lemma 3.4.2]{Kas91}) and they commute with the left multiplication by $f_j$ as follow
$$e_i'f_j=\delta_{ij}+q_i^{-a_{ij}}f_je_i'.$$
These lead to a left $B_q(\g)$-module structure of $U^-$, where $B_q(\g)$ is the quantum Boson algebra with the generators $e_i',f_i (i\in I)$ and the defining relations
$$\begin{aligned}
& e_i'f_j=\delta_{ij}+q_i^{-a_{ij}}f_je_i' \quad\text{for} \ i,j\in I,\\
& \sum_{r+s=1-a_{ij}}(-1)^r
{e_i'}^{(r)}e'_{j}{e_i'}^{(s)}=0 \quad\text{for} \ i\in
I^{\text{re}},j\in I \ \text {and} \ i \neq j, \\
& \sum_{r+s=1-a_{ij}}(-1)^r
{f_i}^{(r)}f_{j}f_i^{(s)}=0 \quad\text{for} \ i\in
I^{\text{re}},j\in I \ \text {and} \ i \neq j, \\
& e_{i}'e_{j}'-e_{j}'e_{i}' = f_{i}f_{j}-f_{j}f_{i} =0 \quad\text{for}\ i,j\in I\ \text{and}\ a_{ij}=0.
\end{aligned}
$$
%Then $U^-$ is a simple $B_q(\g)$-module by the standard argument in \cite{Kas91}.

\vskip 2mm
Define a twisted multiplication on $U^-\otimes U^-$ by
$$(x_1\otimes x_2)(y_1\otimes y_2)=q^{-(|x_2|, |y_1|)}x_1y_1\otimes x_2y_2$$
for homogeneous $x_1,x_2,y_1,y_2$. We have an algebra homomorphism $\rho:U^-\rightarrow U^-\otimes U^-$ given by $\rho{(f_i)}=f_i\otimes 1+1\otimes f_i$ ($i\in I$) with respect to the above algebra structure on $U^-\otimes U^-$. There are two nondegenerate symmetric bilinear forms on $U^-$, namely, Kashiwara's form $( \ , \ )_K$ and Lusztig's form $( \ , \ )_L$, determined by
$$\begin{aligned}
& (1,1)_K=1;\ (e_i'x,y)_K=(x,f_iy)_K,\\
& (1,1)_L=1;\ (e_i'x,y)_L=(1-q_i^2)(x,f_iy)_L
\end{aligned}$$
for $x,y\in U^-$. So that $(f_i,f_i)_K=1$, $(f_i,f_i)_L=(1-q_i^2)^{-1}$ for all $i\in I$.

\vskip 2mm
Let $\mathcal A=\Z[q,q^{-1}]$. The $\mathcal A$-form $_{\mathcal A} U$ is the $\mathcal A$-subalgebra of $U$ generated by the $e_i^{(n)}, f_i^{(n)}$ $(i\in I^\re,n\in \Z_{> 0})$ and $q^h$ $(h\in P^{\vee})$. 

Let $_{\mathcal A} U^-$ (resp. $_{\mathcal A} U^+$) be the $\mathcal A$-subalgebra of $_{\mathcal A} U$ generated by  $f_i^{(n)}$ (resp. $e_i^{(n)}$) and let $^-$ be the $\Q$-algebra involution of $U^-$ given by $\overline{q}=q^{-1}$ and $\overline{f}_i=f_i$. Then the bar-involution $^-$ and the comultiplication $\rho$ restrict to
$$^-: {_{\mathcal A}U}^-\rightarrow {_{\mathcal A} U}^-,\ \ \rho: {_{\mathcal A} U}^-\rightarrow {_{\mathcal A} U}^-\otimes {_{\mathcal A} U}^-.$$
\vskip 2mm

The irreducible highest weight $U_q(\g)$-module $V(\lambda)$ for $\lambda\in P^+$ is given by
$$
\begin{aligned}
 &V(\lambda)\cong U\bigg/(\sum_{i\in I}Ue_i+\sum_{h\in P^\vee}U(q^h-q^{\lambda(h)})+\sum_{i\in I^{\text{re}}}U f_i^{\lambda(h_i)+1}+
\sum_{i\in I^\im \ \text{with}\ \lambda(h_i)=0}Uf_i) \\
 &\phantom{V(\lambda)}\cong U^- \bigg/ %\sum_{(i,l)\in I^{\infty}}U_q(\g)e_{il}+\sum_{h\in P^{\vee}} U_q(\g)(q^h-q^{\lambda(h)})+
(\sum_{i\in I^{\text{re}}}U^- f_i^{\lambda(h_i)+1}+
\sum_{i\in I^\im \ \text{with}\ \lambda(h_i)=0}U^-f_i).
\end{aligned}
$$
The crystal bases of $U^-$ (resp. $V(\lambda)$) constructed in \cite{JKK2005} will be denoted by $ B(\infty)$ (resp. $ B(\lambda)$).
\vskip 2mm
A $U_q(g)$-modules $V$ is said to be integrable  if it satisfies
\begin{itemize}
\item [(i)] $V=\bigoplus_{\mu \in P}V_{\mu}$ with each $\text{dim} V_{\mu} <\infty,$
where $$V_{\mu}=\{v\in V\mid q^{h}=q^{\mu(h)}\ \text{for all}\ h\in P^{\vee}\}.$$
\item [(ii)] $\text{wt}(V):=\{\mu\in P\mid V_{\mu}\neq0\}\subset \bigcup_{i=1}^{s}(\lambda_{i}-Q_{+})$ for finitely many $\lambda_{1},\cdots ,\lambda_{s}\in P$,
\item [(iii)] $e_i$ and $f_{i}$ are locally nilpotent on $V$ for each $i\in I^{\re}$.
\item [(iv)] If $i\in I^{\im}$, then $\mu(h_i)\in\Z_{\geq 0}$ for all $\mu \in \text{wt}(V)$.
\item [(v)] If $i\in I^{\im}$ and $\mu(h_i)=0$, then $f_iV_{\mu}=0$.
\end{itemize}
 If $V$ is an integrable highest weight module with highest weight $\lambda\in P^+$, then
 $V\cong V(\lambda)$.
\vskip 3mm
\subsection{KLR-algebras and the categorification of  ${_{\mathcal A} U}^-$ }\
\vskip 2mm
%Let $\mathcal F$ be the free associative algebra over $\Q(q)$ generated by the symbols $\theta_{i}$ for $i\in I$.  The algebra $\mathcal F$ is $\N[I]$-graded by assigning $|f_i|=i$.
%Let's recall the construction and the main results in \cite{TW2023}.
We  review  the main results in \cite{TW2023} (see \cite{KL2009} for the explanation of the following diagrams).
\vskip 2mm
Given a Borcherds-Cartan datum.  For $\nu=\sum_{i\in I}\nu_ii\in \N[I]$ with $\text{ht}{(\nu)}=\sum_{i\in I}\nu_i=n$, we denote by $\text{Seq}(\nu)$ the set of sequences $\ii=i_1i_2\dots i_n$ in $I$ such that $\nu=i_1+i_2\cdots +i_n$.  The KLR-algebra  $R(\nu)$ is defined to be the $\K$-algebra with generators:

$$1_{\ii}=\genO{i_1}{i_k}{i_n} \quad \text{for} \ \ii=i_1i_2\dots i_n\in \text{Seq}(\nu) \ \text{with} \ \text{deg}(1_{\ii})=0,$$
\vskip 2mm
$$\ \ \ x_{k,\ii}=\genX{i_1}{i_k}{i_n} \quad \text{for} \ \ii\in \text{Seq}(\nu),1\leq k\leq n \ \text{with} \ \text{deg}(x_{k,\ii})=2r_{i_k},$$
\vskip 2mm
$$\ \ \ \ \ \tau_{k,\ii}=\genT{i_1}{i_k}{i_{k+1}}{i_n} \quad \text{for} \ 1\leq k\leq n-1 \ \text{with} \ \text{deg}(\tau_{k,\ii})=-(\alpha_{i_k}, \alpha_{i_{k+1}}).$$
\vskip 4mm
\noindent Subject to the following local relations:

\begin{align}
 \dcross{i}{j} \  =  \  \begin{cases} \quad \quad \quad \quad \quad \  0 & \text{ if } i= j, \\
                    \\
                    \quad  \quad \quad \quad \ \ \dcrossA{i}{j} & \text{ if } i\ne j\ \text{and} \ (\alpha_i,\alpha_j)=0, \\
                    \\
                    \   \dcrossL{-a_{ij}}{i}{j} \ + \ \dcrossR{-a_{ji}}{i}{j}  & \text{ if } i\ne j \ \text{and} \  (\alpha_i,\alpha_j) \ne 0,
                  \end{cases}
\end{align}

\begin{equation}
\begin{aligned}
 \LU{{}}{i}{i} \   -  \ \RD{{}}{i}{i}\ =\ \dcrossA{i}{i} \quad\quad\quad \LD{{}}{i}{i}   \ - \ \RU{{}}{i}{i}\ =\ \dcrossA{i}{i}  \quad \text{ if } i \in I^\re,
\end{aligned}
\end{equation}

\begin{equation}
 \LU{{}}{i}{i} \    =  \ \RD{{}}{i}{i}  \quad\quad\quad\LD{{}}{i}{i}   \  =   \ \RU{{}}{i}{i}  \quad \text{ if } i \in I^\im,
\end{equation}

\begin{equation}
\begin{aligned}
 \LU{{}}{i}{j} \   =  \ \RD{{}}{i}{j}  \quad\quad\quad \LD{{}}{i}{j}   \ =  \ \RU{{}}{i}{j}  \quad \text{ if } i \ne j,
\end{aligned}
\end{equation}

\begin{equation}\label{S1}
\quad \quad\quad\quad \quad\BraidL{i}{j}{i} \  -  \ \BraidR{i}{j}{i} \ =  {\sum_{c=0}^{-a_{ij}-1}} \ \threeDotStrands{i}{j}{i} \ \ \text{ if } i\in I^\re, i\ne j \ \text{and} \ (\alpha_i,\alpha_j) \ne 0,
 \end{equation}
 \begin{equation}
\BraidL{i}{j}{k} \  = \ \BraidR{i}{j}{k} \quad\quad \text{otherwise}.
\end{equation}

\vskip 2mm

\vskip 2mm

For $\ii,\jj\in \Seq(\nu)$, we set $_{\ii}R(\nu)_{\jj}=1_{\ii}R(\nu)1_{\jj}$. Then $_{\ii}R(\nu)_{\jj}$ has a basis $$\{x_{1,\ii}^{u_1}\cdots x_{n,\ii}^{u_n}\cdot\widehat\omega_\jj\mid u_1,\dots,u_n\in \N,\ \omega\in S_n\ \text{such that} \ \omega(\jj)=\ii\},$$ where $\widehat\omega_\jj\in {_{\ii}}\mathcal R(\nu)_{\jj}$ is uniquely determined by a fixed  reduced expression of $\omega$.

\vskip 2mm
Denote by
$$\begin{aligned}
& R(\nu)\text{-}\Mod:\text{the category of finitely generated graded}\ R(\nu) \text{-modules},\\
& R(\nu)\text{-}\fMod:\text{the category of finite-dimensional graded}\ R(\nu) \text{-modules},\\
& R(\nu)\text{-}\pMod:\text{the category of  projective objects in } R(\nu)\text{-}\Mod.
\end{aligned}$$
Let $\B_{\nu}$ be the set of equivalence classes (under isomorphism and degree shifts) of
gr-irreducible $R(\nu)$-modules. Choose one representative $S_b$ from
each equivalence class and denote by $P_b$ the gr-indecomposable projective cover of $S_b$. The Grothendieck group $G_0(R(\nu))$ (resp. $K_0(R(\nu))$) of $R(\nu)$-$\fMod$ (resp. $R(\nu)$-$\pMod$) are free $\Z[q,q^{-1}]$-modules:
  $$G_0(R(\nu))=\bigoplus_{b\in \B_{\nu}}\Z[q,q^{-1}][S_b],\ \ K_0(R(\nu))=\bigoplus_{b\in \B_{\nu}}\Z[q,q^{-1}][P_b].$$
Here, $q[M]=[M\{1\}]$. Let $R=\bigoplus_{\nu\in\N[I]}R(\nu)$, $\B=\bigsqcup_{\nu\in\N[I]}\B_\nu$ and form
$$G_0(R)=\bigoplus_{\nu\in \N[I]}G_0(R(\nu)),\ K_0(R)=\bigoplus_{\nu\in \N[I]}K_0(R(\nu)).$$
\vskip 2mm
The $K_0(R)$ and $G_0(R)$ are dual to each other with respect to the bilinear pairing
$( \ , \ ):K_0(R)\times G_0(R)\rightarrow \Z[q,q^{-1}]$ given by
\begin{equation}\label{KL}([P],[M])=\gdim (P^\psi\otimes_{R(\nu)}M)=\gdim \HOM_{R(\nu)}(\overline{P},M),\end{equation}
where  $\psi$ is the anti-involution of $R(\nu)$ obtained by flipping the diagrams about horizontal axis and it turns a left $R(\nu)$-module into right, $\overline{P}=\HOM(P,R(\nu))^\psi$. There is also a symmetric bilinear form $( \ , \ ):K_0(R)\times K_0(R)\rightarrow \Z ((q))$ defined in the same way.

\vskip 2mm

The $G_0(R)$ and the $K_0(R)$ are equipped with  twisted bialgebras structure induced by the induction and restriction functors defined as follows
\begin{equation*}
\begin{aligned}
& \Ind^{\nu+\nu'}_{\nu,\nu'}:R(\nu)\otimes R(\nu')\text{-}\Mod\rightarrow R(\nu+\nu')\text{-}\Mod,\ M\mapsto R(\nu+\nu')1_{\nu,\nu'}\otimes_{R(\nu)\otimes R(\nu')}M,\\
&\Res^{\nu+\nu'}_{\nu,\nu'}: R(\nu+\nu')\text{-}\Mod\rightarrow R(\nu)\otimes R(\nu')\text{-}\Mod,\ N\mapsto 1_{\nu,\nu'}N,
\end{aligned}
\end{equation*}
where $1_{\nu,\nu'}=1_\nu\otimes 1_{\nu'}\in R(\nu+\nu')$. %For $[M], [N]\in G_0(R)$, we denote the product and the coproduct by
%$$[M]\circ [N]=[ \Ind M\otimes N],\ \ \Delta[M]=[\textstyle \sum \Res M],$$
%and for $[P], [Q]\in K_0(R)$, we write
%$$[P]\cdot [Q]=[ \Ind P\otimes Q],\ \ \rho[P]=[\textstyle \sum \Res P].$$\vskip 2mm
\vskip 2mm

For $\ii\in \Seq(\nu)$, we set $P_{\ii}=R(\nu)1_{\ii}\in R(\nu)\text{-}\pMod$. For $i\in I^\re$ and $n\geq 0$, we set $$P_{i^{(n)}}= R(ni)\psi{(e_{i,n})}\{-\frac{n(n-1)}{2} \cdot r_i\} \cong R(ni)e_{i,n}\{\frac{n(n-1)}{2} \cdot r_i\} ,$$ where $e_{i,n}$ is the primitive idempotent  $x_{1,ni}^{n-1}x_{2,ni}^{n-2}\cdots x_{{n-1},ni}\cdot\tau_{\omega_0}$ of $R(ni)$ with $\omega_0$ be the longest element in $S_n$.
\vskip 2mm
\begin{theorem}\cite{TW2023}\
{\it There is a  $\mathcal A$-bialgebras isomorphism $\gamma: {_{\mathcal A} U}^-\rightarrow K_0(R)$ which sends ${_{\mathcal A} U}^-_{-\nu}$ onto $K_0(R(\nu))$, given by
$$\begin{aligned}& f_i^{(n)}\mapsto [P_{i^{(n)}}]\quad \text{for}\ i\in I^\re,n\geq 0,\\ & f_i\mapsto [P_{i}]\quad \text{for}\ i\in I^\im.\end{aligned}$$
Under this isomorphism, the Lusztig's form $( \ , \ )_L$ corresponds to $( \ , \ )$ on $K_0(R)$, the bar-involution corresponds to $P\mapsto \overline{P}$.
}
\end{theorem}
\vskip 2mm

Let ${_{\mathcal A} U}^{-*}$ be the $\mathcal A$-linear dual of ${_{\mathcal A} U}^-$, which is also a twisted bialgebra with multiplication (resp. comultiplication) dual to the comultiplication (resp. multiplication) of ${_{\mathcal A} U}^-$. The dual map of the bar-involution on ${_{\mathcal A} U}^-$  gives an involution $^-:{_{\mathcal A} U}^{-*}\rightarrow{_{\mathcal A} U}^{-*}$.
\vskip 2mm
%For $M\in R(\nu)$-$\fMod$,  the character of $M$ is defined by $$\Ch M=\sum_{\ii\in \Seq(\nu)}\gdim (1_{\ii}M)\ii\ \in \Z[q,q^{-1}]\Seq(\nu).$$

\begin{theorem}
{\it We identify $K_0(R)^*$ with $G_0(R)$ via   the nondegenerate bilinear pairing in (\ref{KL}). Then the dual map $\gamma^*: G_0(R)\rightarrow {_{\mathcal A} U}^{-*}$ of $\gamma$ is a $\mathcal A$-bialgebras isomorphism.
Under this isomorphism, the bar-involution on ${_{\mathcal A} U}^{-*}$ corresponds to $M\mapsto \HOM_\K(M,\K)^\psi$ of $G_0(R)$.
}
\end{theorem}

\vskip 3mm
\subsection{$R(ni)$-module $V(i^n)$ for $i\in I^\im$}\

\vskip 2mm
For each $i\in I$, the algebra $R(ni)$ for $\nu=ni$ is generated by $x_{1,ni},\dots, x_{n,ni}$ of degree $2r_i$ and $\tau_{1,ni},\dots,\tau_{n-1,ni}$ of degree $-(\alpha_i,\alpha_i)$, subject to the local relations
 $$\dcrossI{i}{i}=0\quad\quad \LU{{}}{i}{i} \   -  \ \RD{{}}{i}{i}\ =\LD{{}}{i}{i}   \ - \ \RU{{}}{i}{i}\ =\ \dcrossA{i}{i}\quad\quad \BraidLI{i}{i}{i}  =  \BraidRI{i}{i}{i} \quad \text{ if } i \in I^\re .$$
 $$\dcrossI{i}{i}=0\quad\quad \LU{{}}{i}{i}    =  \RD{{}}{i}{i}\quad\quad \LD{{}}{i}{i}     =    \RU{{}}{i}{i}\quad\quad \BraidLI{i}{i}{i}  =  \BraidRI{i}{i}{i} \quad \text{ if } i \in I^\im .$$
 \vskip 1mm
We will abbreviate $x_{k,ni}$ (resp. $\tau_{l,ni}$) to $x_k$ (resp. $\tau_l$) since only one sequence is considered. In both cases, $R(ni)$ has a basis $\{x_1^{r_1}\cdots x_n^{r_n}\cdot\tau_\omega\mid \omega\in S_n,r_1,\dots,r_n\geq 0\}$ and $\{\tau_\omega \cdot x_1^{r_1}\cdots x_n^{r_n}\mid \omega\in S_n,r_1,\dots,r_n\geq 0\}$ as an alternative.
We identify the polynomial algebra $P_n=\K[x_1,\dots,x_n]$ with the subalgebra of $R(ni)$ generated by $x_1,\dots,x_n$. Then the center of $R(ni)$ is $P_n^{S_n}$, consisting of  all  symmetric polynomials in $x_1,\dots,x_n$.

\vskip 2mm

Up to isomorphism and degree shifts, each $R{(ni)}$ has a unique gr-irreducible module $V(i^n)$. If $i\in I^\re$,  by the representation theory of the nil-Hecke algebras, $V(i^n)\cong R(ni)\otimes_{P_n}L\{\frac{n(n-1)}{2} \cdot r_i\}$ of graded dimension $[n]_i!$, where $L$ is the one-dimensional trivial module over $P_n$ with each $x_k$ acts by $0$. The gr-projective cover of $V(i^n)$ is $P_{i^{(n)}}$.

Let $\mu=(\mu_1,\dots,\mu_r)$ be any composition of $n$. Define $V({i^\mu}):=V(i^{\mu_1})\otimes \cdots \otimes V(i^{\mu_r})$ to be a gr-irreducible module over the parabolic subalgebra $R(\mu i):=R(\mu_1 i)\otimes \cdots \otimes R(\mu_{r} i)$ of $R(ni)$. Up to  degree shifts, we have
\begin{equation}
\begin{aligned}
& \Ind^{ni}_{\mu i} V(i^\mu)\cong V(i^{n}),\\
& \soc \Res^{ni}_{\mu i} V(i^n)\cong V(i^\mu).
\end{aligned}
\end{equation}
and $\soc \Res^{\ ni}_{(n-1) i} V(i^n)\cong V(i^{n-1})$.

 \vskip 2mm
If $i\in I^\im$, since $R(ni)$ has only trivial idempotents, $V(i^n)$ is the one-dimensional trivial module with the gr-projective cover $R(ni)$.  The following Lemma shows the significant differences with the `real' cases when we use $\Ind$ and $\Res$ to the irreducibles.

\vskip 2mm
For any composition $\mu=(\mu_1,\dots,\mu_r)$ of $n$, denote by $D_{\mu}$ (resp. $D^{-1}_{\mu}$)  the set of minimal length left (resp. right) $S_\mu:=S_{\mu_1}\times \cdots \times S_{\mu_r}$-coset representatives in $S_n$.

\vskip 2mm
\begin{lemma}\label{T}
{\it Let $i\in I^\im$ and $\mu=(\mu_1,\dots,\mu_r)$ be a composition of $n$. Let $M$ be a gr-irreducible $R(\nu)$-module for some $\nu\in \N[I]$. Then
\begin{itemize}
\item[(i)] $\Res^{ni}_{\mu i} V(i^n)\cong V(i^\mu),\ \Res^{\ ni}_{(n-1) i} V(i^n)\cong V(i^{n-1})$.
\item[(ii)] $\Ind^{ni}_{\mu i} V(i^\mu)$ has a unique (graded) maximal submodule $H=\Span\{\tau_u\otimes v\mid u\in D_\mu\backslash \{1\}\}$, where $v$ is a fixed nonzero element in $V(i^\mu)$. The graded head $$\hd \Ind^{ni}_{\mu i} V(i^\mu)\cong V(i^{n}).$$
    Moreover, the $R(\nu)\otimes R(ni)$-module $M\otimes \Ind^{ni}_{\mu i} V(i^\mu)$ has a unique (graded) maximal submodule $M\otimes H$.
\item[(iii)] $\Ind^{ni}_{\mu i} V(i^\mu)$ has a unique (graded) irreducible submodule $V=\Span\{\tau_\lambda\otimes v\}\cong V(i^n)$, where $\lambda$ is the longest shuffle in $D_\mu$. The graded socle $$\soc \Ind^{ni}_{\mu i} V(i^\mu)\cong V(i^{n}).$$
    Moreover, $M\otimes \Ind^{ni}_{\mu i} V(i^\mu)$ has a unique (graded) irreducible submodule $M\otimes V$.
\end{itemize}
}
\begin{proof}
(i) is clear since $V(i^n)$ is $1$-dimensional. Let $N=\Ind^{ni}_{\mu i} V(i^\mu)$.

(ii) For any $u\in D_\mu\backslash \{1\}$, there is no $\omega\in S_n$ such that $\omega u=1$ and $l(\omega u)=l(\omega)+l(u)$. So $H$ is a (graded) maximal submodule of $N$. Note that any nonzero submodule of $\Ind^{ni}_{\mu i} V(i^\mu)$ which contains an element of the form $1\otimes v+\sum_{u\in D_\mu\backslash \{1\}}k_u\tau_u\otimes v$ is equal to $N$. This shows the uniqueness of the maximal. %Let $y_1,\dots,y_k$ be a basis of $M$ consisting of homogenous  elements. Then
%$$W=\Span \{y_i\otimes (\tau_u\otimes v)\mid 1\leq i\leq k,u\in D_\mu\backslash \{1\}\}\cong M\otimes H$$is a graded submodule of $M\otimes N$.
For the second part,  assume $K$ is a nonzero submodule of $M\otimes N$ containing an element
$$y\otimes (1\otimes v)+\sum_{u\in D_\mu\backslash \{1\}}y_u\otimes(\tau_u\otimes v)$$
for some $y, y_u\in M$ and $y\neq 0$. Since the ungraded  $R(\nu)$-module $\underline M$ is simple, we have $M=R(\nu)y$. Thus $y\otimes (1\otimes v)\in K$ and so $K=M\otimes N$, which implies $M\otimes H$ is the unique maximal submodule of $M\otimes N$.

 (iii) Let $L$ be a nonzero submodule of $N$ and $m=\sum_{u\in D_\mu}k_u\tau_u\otimes v\in L$, $m\neq 0$. Choose an element $u$ of  minimal length  such that $k_u\neq 0$. Let $\lambda=\omega u$. Then $l(\lambda)=l(\omega)+l(u)$. If $u'$ satisfies
 $$k_{u'}\neq 0,\ \omega u'=z\in D_{\mu} \ \text{and}\ l(z)=l(\omega)+l(u'), $$
 then $l(z)\leq l(\lambda), l(u')\geq l(u)$ implies $z=\lambda$ and $u=u'$. Therefore $\tau_\omega\cdot m=k_u\tau_\lambda\otimes v\in L$ and so $V\subseteq L$. Note that $V$ is a submodule of $N$ since for any $\omega$ such that $\omega\lambda\in D_\mu$ and $l(\omega\lambda)=l(\omega)+l(\lambda)$, we must have $\omega=1$. This also shows that $V$ is a trivial module, which is isomorphic to $V(i^n)$.

 For the second part, let $K$ be a nonzero submodule of $M\otimes N$. Choose a nonzero $$z=\sum_{u\in D_\mu}y_u\otimes (\tau_u\otimes v)\in K$$ for some $ y_u\in M$. Let $u$ is of minimal length  such that $y_u\neq 0$ and let $\lambda=\omega u$. The same argument as above  leads to
 $$1_\nu\otimes \tau_\omega\cdot z=y_u\otimes (\tau_\lambda\otimes v)\in K.$$
 It follows that $M\otimes V\subseteq K$.
\end{proof}
\end{lemma}

\vskip 6mm

\section{\textbf{Crystal operators $\widetilde{f}_i$ and $\widetilde{e}_i$ } }
\vskip 2mm

For $\nu,\nu'\in\N[I]$, we define the co-induction functor  by
$$\Indd^{\nu+\nu'}_{\nu,\nu'}:R(\nu)\otimes R(\nu')\text{-}\Mod\rightarrow R(\nu+\nu')\text{-}\Mod,\ M\mapsto \HOM_{R(\nu)\otimes R(\nu')}(1_{\nu,\nu'}R(\nu+\nu'),M),$$
where the $R(\nu+\nu')$-module structure of $\Indd^{\nu+\nu'}_{\nu,\nu'}M$ is given by
$$(z\cdot f)(x)=f(xz)\ \ \text{for}\ z\in R(\nu+\nu'), x\in 1_{\nu,\nu'}R(\nu+\nu').$$

The co-induction functor  is right adjoint to the restriction
\begin{equation}\label{adjoint}\HOM_{\nu+\nu'}(M,\Indd_{\nu,\nu'}N)\cong \HOM_{\nu,\nu'}(\Res_{\nu,\nu'} M,N)\end{equation}
for $M\in R(\nu)\otimes R(\nu')\text{-}\Mod, N\in R(\nu+\nu')\text{-}\Mod$. It follows from \cite[Thoerem 2.2]{LV2009} (the proof there also applies in our case) that for $M\in R(\nu)$-$\fMod$ and $N\in R(\nu')$-$\fMod$, we have
\begin{equation}\label{rev}\Ind_{\nu,\nu'}M\otimes N\cong \Indd_{\nu',\nu}N\otimes M \{-\nu\cdot\nu'\}\end{equation}

\vskip 2mm
For $i\in I$ and $n\geq 0$, define the functor  $$\Delta_{i^n}: R(\nu)\text{-}\Mod\rightarrow R(\nu-ni)\otimes R(ni)\text{-}\Mod,\ M\mapsto (1_{\nu-ni}\otimes 1_{ni})M.$$
So we have for $M\in R(\nu)$-$\Mod$ and $N\in R(\nu-ni)$-$\Mod$,
\begin{equation}\label{adj}\HOM_{\nu}(M,\Indd_{\nu-ni,ni}N)\cong \HOM_{\nu-ni,ni}(\Delta_{i^n} M,N).\end{equation}
Let $\varepsilon_i(M)=\Max\{n\geq 0\mid \Delta_{i^n}M\neq 0\}$ be the number of the largest $i$-tail in sequence $\kk$ such that $1_\kk M\neq 0$.
\vskip 2mm

\vskip 2mm

The following three lemmas have been given in \cite[Section 3.2]{KL2009}.

\vskip 2mm

\begin{lemma}\label{L1}
{\it Let $i\in I$ and  $M\in R(\nu)$-$\fMod$ be a gr-irreducible module. If $N\otimes V(i^n)$ is a gr-irreducible submodule of $\Delta_{i^n}M$ for some $0\leq n\leq \varepsilon_{i}(M)$, then $\varepsilon_i(N)=\varepsilon_i(M)-n$.
}
\end{lemma}
\vskip 2mm
\begin{lemma}\label{L2}
{\it Let $i\in I$ and $K\in R(\nu)$-$\fMod$ be a gr-irreducible module with $\varepsilon_i(K)=0$. Set $M=\Ind_{\nu,ni}K\otimes V(i^n)$. Then
\begin{itemize}
\item[(i)] $\Delta_{i^n}M\cong K\otimes V(i^n)$,
\item[(ii)] $\hd M$ is gr-irreducible with $\varepsilon_i(\hd M)=n$,
\item[(iii)] all other composition factors $L$ of $M$ have $\varepsilon_i(L)<n$.
\end{itemize}
}
\end{lemma}

\vskip 2mm
\begin{lemma}\label{LL}
{\it Let $i\in I$ and  $M\in R(\nu)$-$\fMod$ be a gr-irreducible module with $\varepsilon_i(M)=a$. Then $\Delta_{i^a}M$ is isomorphic to $K\otimes V(i^a)$ for some gr-irreducible $K\in R(\nu-ai)$-$\fMod$ with $\varepsilon_i(K)=0$. Furthermore, $M\simeq \hd\Ind_{\nu-ai,ai}K\otimes V(i^a)$ in this case.
}
\end{lemma}
\vskip 2mm

The following lemma is also well known. The proof is dual to Lemma \ref{L2}.
\vskip 2mm
\begin{lemma}\label{L3}
{\it Let $i\in I$ and $K\in R(\nu)$-$\fMod$ be a gr-irreducible module with $\varepsilon_i(K)=0$. Set $M=\Indd_{\nu,ni}K\otimes V(i^n)$. Then
\begin{itemize}
\item[(i)] $\Delta_{i^n}M\cong K\otimes V(i^n)$,
\item[(ii)] $\soc M$ is gr-irreducible with $\varepsilon_i(\soc M)=n$,
\item[(iii)] all other composition factors $L$ of $M$ have $\varepsilon_i(L)<n$.
\end{itemize}
}
\begin{proof}

(i) By (\ref{rev}), $\Delta_{i^n}M\cong \Delta_{i^n}\Ind_{ni,\nu} V(i^n)\otimes K\{\nu\cdot ni\}$. Then by the Mackey Theorem \cite[Proposition 3.4]{TW2023} and $\varepsilon_i(K)=0$, $\Delta_{i^n} M$ is isomorphic to
$$\Ind^{\nu,ni}_{0,\nu,ni,0}{}^\diamond(\Res^{ni,\nu}_{0,ni,\nu,0}V(i^n)\otimes K)\cong K\otimes V(i^n).$$

(ii) For each nonzero graded submodule $N$ of $M$, by  (\ref{adj}) and the irreducibility of $K\otimes V(i^n)$, we have $\Delta_{i^n} N\twoheadrightarrow K\otimes V(i^n)$.  Suppose we are given a  decomposition
 $$\soc M=N_1\oplus N_2 \oplus \cdots \oplus N_s,$$ such that each $N_k$ is gr-irreducible. Then $K\otimes V(i^n)$ is a quotient of each  $\Delta_{i^n}(N_k)$ and  $\Delta_{i^n}(\soc M)$. But they are nonzero submodules of $\Delta_{i^n}M$, which is irreducible by (i). So $\Delta_{i^n}(\soc M)\cong \Delta_{i^n}(N_k) \cong K\otimes V(i^n)$ and so $\soc M$ must be gr-irreducible.

(iii) Our assertion follows from the irreducibility of $\soc M$ and the exactness of $\Delta_{i^n}$.
\end{proof}
\end{lemma}

\vskip 2mm
\begin{corollary}\label{C1}
{\it Let $i\in I$ and $K\in R(\nu)$-$\fMod$ be  gr-irreducible  with $\varepsilon_i(K)=0$. Then
$$\soc\Indd_{\nu,ni}K\otimes V(i^n)\cong \hd\Ind_{\nu,ni}K\otimes V(i^n).$$
}\begin{proof}
By Lemma \ref{L2}, $\Delta_{i^n}\hd\Ind_{\nu,ni}K\otimes V(i^n)\cong K\otimes V(i^n)$. So by (\ref{adjoint}), we have
$$\hd\Ind_{\nu,ni}K\otimes V(i^n)\hookrightarrow \Indd_{\nu,ni}K\otimes V(i^n).$$
Now the result follows from the previous lemma.
\end{proof}\end{corollary}
 \vskip 2mm
 The following lemma has appeared in \cite[Lemma 3.9]{KOP2012}. Our proof is different and  includes the case where $a_{ii}=0$.

\vskip 2mm
\begin{lemma}\label{L4}
{\it Let $i\in I^\im$ and $K\in R(\nu)$-$\fMod$ be a gr-irreducible module with $\varepsilon_i(K)=0$. Let $\mu=(\mu_1,\dots,\mu_r)$ be a composition of $n$. Set $$M=\Ind^{\nu+ni}_{\nu,\mu i}K\otimes V(i^\mu).$$ Then $\hd M$ is gr-irreducible with $\varepsilon_i(\hd M)=n$. Moreover,
$\hd M\cong \hd\Ind^{\nu+ni}_{\nu,n i} K\otimes V(i^n)$.
}
\begin{proof}

By the  Mackey Theorem  and $\varepsilon_i(K)=0$, $\Delta_{i^n} M\cong K\otimes \Ind^{ni}_{\mu i}V(i^\mu)$. For each  nonzero quotient $Q$ of $M$, $\Delta_{i^n}Q$ is a nonzero quotient of $\Delta_{i^n}M$ (we have a nonzero homomorphism from $K\otimes \Ind^{ni}_{\mu i}V(i^\mu)$ to $\Delta_{i^n} Q$ by Frobenius reciprocity). So $\hd\Delta_{i^n}Q$ is a nonzero quotient of $\hd\Delta_{i^n}M$. By Lemma \ref{T}(ii), $\hd\Delta_{i^n}M\cong K\otimes V(i^n)$. Thus,
$$\hd\Delta_{i^n}Q\cong \hd\Delta_{i^n}M \cong K\otimes V(i^n).$$

Assume we have a decomposition
 $$\hd M=M/M_1\oplus M/M_2 \oplus \cdots \oplus M/M_t,$$ such that each $M/M_k$ is gr-irreducible. Then $$\hd(\Delta_{i^n}\hd M)\cong \bigoplus^{t}_{k=1}\hd(\Delta_{i^n} M/M_k).$$
 But $\hd(\Delta_{i^n}\hd M)$ and each $\hd(\Delta_{i^n} M/M_k)$ are all isomorphic to $K\otimes V(i^n)$, which implies $\hd M$ must be gr-irreducible. So $M$ contains a unique graded maximal submodule $J^{gr}(M)$. Recall from Lemma \ref{T}(ii)  that $H$ denotes the unique (graded) maximal submodule of $N= \Ind^{ni}_{\mu i} V(i^\mu)$ and $N/H\cong V(i^n)$. So we have $\Ind_{\nu,ni}^{\nu+ni} K\otimes H\subseteq J^{gr}(M)$ and a surjection $\Ind_{\nu,ni}^{\nu+ni} K\otimes V(i^n)\twoheadrightarrow \hd M$. By the exactness of $\Delta_{i^n}$ and Lemma \ref{L2}(i), $K\otimes V(i^n)\twoheadrightarrow \Delta_{i^n}\hd M$. Therefore, $ \Delta_{i^n}\hd M\cong K\otimes V(i^n)$ and by Lemma \ref{LL}, $\hd M\cong \hd\Ind^{\nu+ni}_{\nu,n i} K\otimes V(i^n)$.
 \end{proof}
\end{lemma}
\vskip 2mm

By the duality, we have:

\vskip 2mm
\begin{lemma}\label{L5}
{\it Let $i\in I^\im$ and $K\in R(\nu)$-$\fMod$ be a gr-irreducible module with $\varepsilon_i(K)=0$. Let $\mu=(\mu_1,\dots,\mu_r)$ be a composition of $n$. Set $$M={\Indd}^{\nu+ni}_{\nu,\mu i}K\otimes V(i^\mu).$$ Then $\soc M$ is gr-irreducible with $\varepsilon_i(\soc M)=n$. Moreover,
$\soc M\cong \soc\Indd^{\nu+ni}_{\nu,n i} K\otimes V(i^n)$.
}
\begin{proof}
Note that $M\cong \Ind_{ni,\nu}(\Ind V(i^{\mu'}))\otimes K$, where $\mu'=(\mu_r,\dots,\mu_1)$.

By the  Mackey Theorem  and $\varepsilon_i(K)=0$, $\Delta_{i^n} M\cong K\otimes \Ind V(i^{\mu'})$.  For each nonzero graded submodule $N$ of $M$, $\Delta_{i^n}N$ is a nonzero submodule of $\Delta_{i^n}M$ (we have a nonzero homomorphism from $\Delta_{i^n}N$ to $K\otimes \Ind V(i^{\mu'})$  by (\ref{adjoint})). So
$\soc\Delta_{i^n}N \hookrightarrow \soc\Delta_{i^n}M$ is nonzero.
By Lemma \ref{T}(iii), $\soc\Delta_{i^n}M\cong K\otimes V(i^n)$. Thus,
$$\soc\Delta_{i^n}N\cong \soc\Delta_{i^n}M \cong K\otimes V(i^n).$$

Assume we have a decomposition
 $$\soc M=N_1\oplus N_2 \oplus \cdots \oplus N_s,$$ such that each $N_k$ is gr-irreducible. Then $$\soc(\Delta_{i^n}\soc M)\cong \bigoplus^{s}_{k=1}\soc(\Delta_{i^n} N_k).$$ Hence $\soc M$ must be gr-irreducible. Recall that  the unique (graded) irreducible submodule of $\Ind^{ni}_{\mu' i} V(i^{\mu'})$ is isomorphic to $V(i^n)$. So we have $\soc M\hookrightarrow \Indd^{\nu+ni}_{\nu,ni}K\otimes V(i^n)$. By the exactness of $\Delta_{i^n}$ and Lemma \ref{L3}(i), $\Delta_{i^n}\soc M\hookrightarrow K\otimes V(i^n)$. Therefore $ \Delta_{i^n}\soc M\cong K\otimes V(i^n)$ and then by Lemma \ref{LL} and Corollary \ref{C1}, $\soc M\cong \hd\Ind_{\nu,n i} K\otimes V(i^n)\cong \soc\Indd_{\nu,ni}K\otimes V(i^n)$.
\end{proof}
\end{lemma}

\vskip 2mm
\begin{proposition}\label{P1}
{\it Let $i\in I$ and $N\in R(\nu)$-$\fMod$ be gr-irreducible with $\varepsilon_i(N)=a$. Assume $\Delta_{i^a}N\cong K\otimes V(i^a)$ for  gr-irreducible $K\in R(\nu-ai)$-$\fMod$ with $\varepsilon_i(K)=0$. Let $M=\Ind^{\nu+ni}_{\nu,ni}N\otimes V(i^n)$. Then
 $$\hd M\cong \hd \Ind K\otimes V(i^{a+n}),$$
which is gr-irreducible with $\varepsilon_i(\hd M)=\varepsilon_i(N)+n$. Moreover, if $i\in I^\re$, all other composition factors of $M$ with $\varepsilon_i< \varepsilon_i(N)+n$
}
\begin{proof}
The case when $i\in I^\re$ has been shown in \cite[Section 3.2]{KL2009}. Suppose $i\in I^\im$. Since $N$ is a quotient of $\Ind K\otimes V(i^a)$, so by the the exactness and the  transitivity of induction, $M$ is a quotient of
 $\Ind K\otimes(\Ind V(i^a)\otimes V(i^n))$. Hence, $\hd M$ is a quotient of
 $\hd \Ind K\otimes(\Ind V(i^a)\otimes V(i^n))$.  Now the result
 follows from Lemma \ref{L4}.
\end{proof}
\end{proposition}

\vskip 2mm
\begin{remark}
If $i\in I^\im$, $M=\Ind V(i)\otimes V(i)$ has a composition series
$$M\supsetneq J^{gr}(M)=\Span\{\tau_1\otimes v\}\supsetneq 0$$
with all composition factors isomorphic to $V(i^2)$. So the last part of Proposition \ref{P1} is not true for the imaginary cases.
\end{remark}

\vskip 2mm
\begin{proposition}\label{P2}
{\it Let $i\in I$ and $N\in R(\nu)$-$\fMod$ be gr-irreducible with $\varepsilon_i(N)=a$. Assume $\Delta_{i^a}N\cong K\otimes V(i^a)$ for  gr-irreducible $K\in R(\nu-ai)$-$\fMod$ with $\varepsilon_i(K)=0$. Let $M=\Indd^{\nu+ni}_{\nu,ni}N\otimes V(i^n)$. Then
 $$\soc M\cong \soc \Indd K\otimes V(i^{a+n})\ (\ \cong \hd \Ind K\otimes V(i^{a+n})\cong \hd \Ind N\otimes V(i^{n})\ ),$$
which is gr-irreducible with $\varepsilon_i(\soc M)=\varepsilon_i(N)+n$. Moreover, if $i\in I^\re$, all other composition factors of $M$ with $\varepsilon_i< \varepsilon_i(N)+n$
}
\begin{proof}
Dualizing the proof of Proposition \ref{P1} and using Lemma \ref{L5} for $i\in I^\im$.
\end{proof}
\end{proposition}

\vskip 2mm
Assume $i\in I$ and $M\in R(\nu)$-$\fMod$ be gr-irreducible. Assume further that $\varepsilon_i(M)=a$ and  $\Delta_{i^a}M\cong K\otimes V(i^a)$ for some gr-irreducible $K\in R(\nu-ai)$-$\fMod$ with $\varepsilon_i(K)=0$. By Proposition \ref{P1}, we can define
$$\widetilde{f}_{i}M:=\hd\Ind^{\nu+i}_{\nu,i}M\otimes V(i)\cong \hd\Ind^{\nu+i}_{\nu-ai,(a+1)i} K\otimes V(i^{a+1}).$$
Furthermore, we deduce from Proposition \ref{P1} and  \ref{P2} that, for any $n\geq 1$,
\begin{equation}\begin{aligned}
& {\widetilde{f}_{i}}^nM\cong\hd\Ind M\otimes V(i^n)\cong \hd\Ind K\otimes V(i^{a+n})\\
& \phantom{{\widetilde{f}_{i}}^nM}\cong \soc \Indd M\otimes V(i^{n})\cong \soc \Indd K\otimes V(i^{a+n}).
\end{aligned}
\end{equation}
In particular, $M\cong {\widetilde{f}_{i}}^aK$.
For the same $M$, another crystal operator $\widetilde{e}_i$ is defined to be $$\widetilde{e}_iM:=\soc (e_i M).$$
 where $e_i=\Res_{\nu-i}^{\nu-i,i}\circ \Delta_i:R(\nu)\text{-}\fMod\rightarrow R(\nu-i)\text{-}\fMod$. Then $\widetilde{e}_iM$ is either gr-irreducible with $\varepsilon_i(\widetilde{e}_iM)=a-1$ when $a>0$, or equal to $0$ when $a=0$. It was proved in \cite[Chapter 4]{TW2023} that, for any $1\leq n\leq a$,
 $$\widetilde{e}_{i}^nM\cong \hd\Ind K\otimes V(i^{a-n}).$$
 So we see that $K\cong \widetilde{e}_{i}^aM$ and
\begin{equation}\widetilde{e}_{i}\widetilde{f}_{i}M\cong M;\ \widetilde{f}_{i}\widetilde{e}_{i}M\cong M \ \text{if} \ \varepsilon_i(M)>0.\end{equation}

\vskip 2mm

Using the fact proved in \cite[Chapter 4]{TW2023} that
 \begin{equation}\label{aaa}
 \soc(\Delta_{i^n}M)\cong\widetilde{e}_{i}^nM\otimes V(i^n),
 \end{equation}
we have the following result.

 \vskip 2mm
\begin{proposition}
{\it Let   $M\in R(\nu)$-$\fMod$ be  gr-irreducible  and let $n\geq 0$. If $i\in I^\re$, then
 $$\soc(e_i^nM)\cong (\widetilde{e}_{i}^nM)^{[n]_i!}.$$
 If $i\in I^\im$, then $\soc(e_i^nM)\cong \widetilde{e}_{i}^nM$.
}
\begin{proof}
The case of $i\in I^\re$ has been proved in \cite[Lemma 3.14]{KL2009}. We can assume $i\in I^\im$.

 Assume  $\varepsilon_i(M)=a$ and  $\Delta_{i^a}M\cong K\otimes V(i^a)$ for a gr-irreducible $K\in R(\nu-ai)$-$\fMod$ with $\varepsilon_i(K)=0$. Then $M\cong\hd\Ind K\otimes V(i^a)$. By Mackey's Theorem, we have $$\begin{aligned}& \Delta_{i^n} \Ind K\otimes V(i^a) \cong\Ind^{\nu-ni,ni}_{\nu-ai,(a-n)i,ni}K\otimes \Delta_{i^n}V(i^a) \\
&\phantom{\Delta_{i^n} \Ind K\otimes V(i^a) }\cong\left(\Ind^{\nu-ni}_{\nu-ai,(a-n)i}K\otimes V(i^{a-n})\right)\otimes V(i^n).\end{aligned}$$
Hence $1_{\nu-ni}\otimes R(ni)\cdot\Delta_{i^n}M=0$. Let $L$ be a gr-irreducible submodule of $e_i^nM$, then $L$ is an $R(\nu-ni)\otimes R(ni)$-submodule of $\Delta_{i^n}M$, which is isomorphic to $L\otimes V(i^n)$. Since $\soc(\Delta_{i^n}M)$ is gr-irreducible, we must have $L=\soc(e_i^nM)$.
Apply (\ref{aaa}), we  have an isomorphism of $R(\nu-ni)\otimes R(ni)$-modules
 $$L\otimes V(i^n)\cong L=\soc(\Delta_{i^n}M)\cong\widetilde{e}_{i}^nM\otimes V(i^n).$$
 Thus $L\cong \widetilde{e}_{i}^nM$ as $R(\nu-ai)$-modules. This completes the proof.
\end{proof}
\end{proposition}
\vskip 6mm
%where $P=\Ind K\otimes V(i^a)$ and $J=J^{gr}(P)$. Let $m=\text{ht}{(\nu-ai)}$. Then we have
%$$P=\Span\{\widehat\omega_\jj\otimes (K\otimes V(i^a))\mid \omega\in D_{(m,a)},\ \jj=\textit{\textbf{l}} i^a \ \text{with}\ \textit{\textbf{l}} \in \Seq(\nu-ai)\ \text{such that}\   l_m\neq i \},$$
%and further $$\Delta_{i^n}P=\Span\{\widehat\omega_\jj\otimes (K\otimes V(i^a))\mid \omega\in D_{(m,a)}\cap S_{m+a-n},\ \jj=\textit{\textbf{l}} i^a \ \text{with}\   l_m\neq i \}.$$
\section{\textbf{Categorification of crystals $B(\infty)$ and $ B(\lambda)$ } }

%In this section we closely follow the frameworks given in \cite{LV2009} and \cite{KOP2012}.
\vskip 2mm
By the symmetry of the KLR-algebras, for a gr-irreducible $M\in R(\nu)$-$\fMod$, we similarly define
$$
\begin{aligned}
& \Delta_i^\vee (M)=\Res^{\nu}_{i,\nu-i}M, \ \ e_i^\vee (M)=\Res_{\nu-i}^{i,\nu-i}\Delta_i^\vee (M),\\
& \widetilde{e}_{i}^\vee(M)=\soc(e_i^\vee M), \ \ {\widetilde{f}_{i}}^\vee(M)=\hd \Ind_{i,\nu}^{\nu+i} V(i)\otimes M,\\
&\varepsilon_i^\vee(M)=\Max\{m\geq 0\mid ({\widetilde{e}_{i}^\vee})^m M\neq 0\}.
\end{aligned}
$$
The following lemma is an analogue of Theorem 5.5.1 in \cite{K2005}.
\vskip 2mm
\begin{lemma}\label{L6}
{\it Let  $M\in R(\nu)$-$\fMod$ be  gr-irreducible with $\varepsilon_i^\vee(M)=a$. If $i\in I^\re$,
$$[e_i^\vee M]=[a]_i[\widetilde{e}_{i}^\vee M]+\textstyle\sum c_r[N_r]$$
for some  gr-irreducible $N_r\in R(\nu-i)$-$\fMod$ with $\varepsilon_i^\vee(N_r)<a$. Moreover, $a$ is the maximal size of a Jordan block of $x_{1,\ii}$ on $M$ (with
 eigenvalue $0$) for some $\ii\in \Seq(\nu)$ with $i_1=i$.

If $i\in I^\im$, we have $[e_i^\vee M]=[\widetilde{e}_{i}^\vee M]+\textstyle\sum c_r[N_r]$, and $x_{1,\ii}M=0$ for any $\ii\in \Seq(\nu)$ with $i_1=i$.
}
\end{lemma}
\vskip 2mm
For $i\in I$ and $u\in U^-$, we set $l_i(u)=\text{max}\{m\geq 0\mid (e_i')^{m}u\neq 0\}$ and  $$U_i^{-<k}=\{u\in U^-\mid l_i(u)<k\}.$$
\vskip 2mm
\begin{definition}{
A perfect basis of $U^-$ is a basis $B$ consisting of weight
 vectors, such that for any $i\in I$,
\begin{itemize}
\item [(1)] if $b\in B$ with $e_i'(b)\neq 0$, then there is a unique $\mathtt e_i(b)\in B$ such that $$e_i'(b)=c\cdot \mathtt e_i(b)+U^{-<l_i(b)-1}_i \ \text{for some} \ c
\in \Q(q)\backslash \{0\},$$
\item [(2)] if $\mathtt e_i(b)=\mathtt e_i(b')$, then $b=b'$.
\end{itemize}
}\end{definition}
\vskip 2mm
The existence of the perfect basis follows from  the existence of the `upper global basis' of $U^-$ (see \cite{KOP2012}). If we set $\text{wt}(b)=\mu\in Q_-$ for $b\in B_{\mu}$ and
$$\begin{aligned}
& \mathtt f_i(b)=\begin{cases} b' \ \ \text{if}\ \mathtt e_i(b')=b, \\ 0\quad \text{otherwise} , \quad \end{cases} \quad \varepsilon_i(b)=\begin{cases} l_i(b)\ \ \text{if}\ i\in I^\re, \\ 0\qquad \text{if}\ i\in I^\im,\end{cases}\\
& \varphi_i(b)=\varepsilon_i(b)+ \text{wt}_i(b),\ \ \text{where}\ \text{wt}_i(b)=\text{wt}(b)(h_i)
\end{aligned}$$
Then  $(B, \text{wt}, \mathtt e_i,\mathtt f_i, \varepsilon_i,\varphi_i)$ is an abstract crystal defined in \cite{JKKS2007}. Moreover, it was shown in \cite{KOP2012} that the abstract crystal arising from  any perfect basis of $U^-$ is isomorphic to $B(\infty)$. We next show that the set of isomorphic classes $\mathbb B$ of gr-irreducible $R$-modules is a perfect basis of $U^-$, so that $\mathbb B\cong  B(\infty)$ as abstract crystals.
\vskip 2mm

For $[M]\in G_0(R)$, we define $\mathbf f_i[M]=[\Ind V(i)\otimes M]$. Then by  the Mackey Theorem, we have the following equation in $G_0(R)$:
$$e_i^\vee \mathbf f_j[M]=\delta_{ij}+q_i^{-a_{ij}} \mathbf f_j e_i^\vee[M].$$
Let $[M]\in G_0(R), [P]\in K_0(R)$. We have by the property of  pairing (\ref{KL}) that
\begin{equation}\label{dual}
([P],e_i^\vee[M])=(f_i[P],[M]), \ \ ([P],\mathbf f_i[M])=(e_i'[P],[M]).
\end{equation}
It follows that $G_0(R)_{\Q(q)}$ is a $B_q(\g)$-module with the action of $e_i'$ by $e_i^\vee$ and $f_i$ by $\mathbf f_i$. Note we have a linear isomorphism
$$\psi:G_0(R)_{\Q(q)}\rightarrow U^-,\ \ [M]\mapsto u_M$$
such that, for any $[P]\in U^-$, $([P],[M])=(u_M,[P])_K$. Then by (\ref{dual}),
$$(u_{e_i^\vee M}, [P])_K=([P], e_i^\vee [M])=(f_i[P],[M])=(u_M,f_i[P])_K=(e_i'u_{M},[P])_K,$$
$$(u_{\mathbf f_i M}, [P])_K=([P], \mathbf f_i [M])=(e_i'[P],[M])=(u_M,e_i'[P])_K=(f_i u_{M},[P])_K,$$
thus $u_{e_i^\vee M}= e_i'u_{M}$ and $u_{\mathbf f_i M}=f_i u_{M}$, which implies $\psi$ is $B_q(\g)$-linear. Then by Lemma \ref{L6}, $\mathbb B$ is a perfect basis of $B_q(\g)$-module $G_0(R)_{\Q(q)}$ with the data
$$\mathtt e_i([M])=\widetilde{e}_{i}^\vee [M],\ \ l_i([M])=\varepsilon_i^\vee(M).$$

\vskip 2mm
For $\lambda\in P^+$ and $i\in I$, we set $\lambda_i=\lambda(h_i)\geq 0$. The cyclotomic KLR-algebra $R^{\lambda}(\nu)$ is defined to be the quotient of $R(\nu)$ by  $\mathcal J^\lambda_\nu$, which is the two sided ideal of $R(\nu)$  generated by $x_{1,\ii}^{\lambda_{i_1}}$ for all $\ii\in I$. Let  $R^\lambda=\bigoplus_{\nu}R^{\lambda}(\nu)$.

\vskip 2mm

We have an obvious corollary of Lemma \ref{L6}.
\vskip 2mm
\begin{corollary}\label{CC}
{\it Let  $\Lambda\in P^+$, $M\in R(\nu)$-$\fMod$ be  gr-irreducible. Then $\mathcal J^\lambda_\nu M=0$  if and only if $\varepsilon_i^\vee(M)\leq \lambda_i$ for all $i\in I^\re$ and $\varepsilon_i^\vee(M)=0$ for all $i\in I^\im$ with $\lambda_i=0$.
}
\end{corollary}
\vskip 2mm
Define the  projection functor by
$$\pr_\lambda : R(\nu)\text{-}\fMod\rightarrow R^{\lambda}(\nu)\text{-}\fMod,\ \ M\mapsto M/\mathcal J^\lambda_\nu M,
$$
and  the inflation  functor $\infl_\lambda :R^{\lambda}(\nu)\text{-}\fMod\rightarrow R(\nu)\text{-}\fMod$ by regarding an $R^{\lambda}(\nu)$-module $\mathcal M$ as an $R(\nu)$-module.

Let $\mathbb B^\lambda$ be the set of equivalence classes (under isomorphism and degree shifts) of gr-irreducible $R^{\lambda}$-modules. Then $\mathbb B^\lambda$ is a basis of  the Grothendieck group $G_0(R^\lambda)$ of $R^\lambda$-$\fMod$. For $\mathcal M\in \mathbb B^\lambda$, we set $\text{wt}^\lambda (\mathcal M)=\lambda-\nu$ and
$$\begin{aligned}
&  e_i^\lambda:R^{\lambda}(\nu)\text{-}\fMod\rightarrow R^{\lambda}(\nu-i)\text{-}\fMod,\ \ \mathcal M\mapsto\pr_\lambda\circ e_i\circ \infl_\lambda \mathcal M,\\
& \widetilde{e}_{i}^\lambda:\mathbb B^\lambda\rightarrow \mathbb B^\lambda\sqcup\{0\},\ \ \mathcal M\mapsto\pr_\lambda\circ \widetilde{e}_{i}\circ \infl_\lambda \mathcal M,  \\ &  \widetilde{f}_{i}^\lambda:\mathbb B^\lambda\rightarrow \mathbb B^\lambda\sqcup\{0\},\ \ \mathcal M\mapsto\pr_\lambda\circ \widetilde{f}_{i}\circ \infl_\lambda \mathcal M,\\
& \varepsilon_i^\lambda(\mathcal M)=\begin{cases}\ \varepsilon_i(\infl_\lambda\mathcal M)  \quad\text{if}\ i\in I^\re\\ \ 0 \qquad \text{if}\ i\in I^\im\end{cases}\\
& \varphi_i^\lambda(\mathcal M)=\begin{cases}\text{max}\{m\geq0\mid \pr_\lambda\circ \widetilde{f}_{i}^m\circ \infl_\lambda \mathcal M\neq 0\}\quad \text{if}\ i\in I^\re\\  \text{wt}_i^\lambda(\mathcal M) \qquad \text{if}\ i\in I^\im .\end{cases}
\end{aligned}$$
It was proved in \cite{LV2009} that if $i\in I^\re$, then $0\leq\varphi_i^\lambda(\mathcal M)<\infty$ and $\varphi_i^\lambda(\mathcal M)=\varepsilon_i^\lambda(\mathcal M)+\text{wt}_i^\lambda(\mathcal M)$. Therefore, $(\mathbb B^\lambda, \text{wt}^\lambda, \widetilde{e}_{i}^\lambda,\widetilde{f}_{i}^\lambda, \varepsilon_i^\lambda,\varphi_i^\lambda)$ is an abstract crystal. The same proof as in \cite[Lemma 5.13]{KOP2012} shows that if $i\in I^\im$,
$\text{wt}_i^\lambda(\mathcal M)=0$ if and only if $\widetilde{f}_{i}^\lambda(\mathcal M)=0$. Then we see by a similar argument in \cite[Theorem 5.14]{KOP2012} that $\mathbb B^\lambda\cong B(\lambda)$ as abstract crystals.

\vskip 2mm

Let $\lambda\in P^+$. Denote the ($\mathcal A$-linear) graded duals of $G_0(R)$ and $G_0(R^\lambda)$ by
$$G_0^*(R)=\bigoplus_\nu G_0(R(\nu))^*,\ \ G_0^*(R^\lambda)=\bigoplus_\nu G_0(R^\lambda(\nu))^*.$$
Note that $G_0(R)$ and $G_0(R^\lambda)$ are natural ${_{\mathcal A} U}^+$-modules, where  $e_i\in {_{\mathcal A} U}^+$ acts on $G_0(R^\lambda)$ by $e_i^\lambda$. Let $*$ be the $\mathcal A$-linear anti-involution such that $e_i^*=e_i$ for all $i\in I$. Then $G_0^*(R)$, similarly for $G_0^*(R^\lambda)$, is a left $_{\mathcal A} U^+$-module with the action given by
\begin{equation}\label{action}
(z\cdot f)([M])=f(z^*[M])\quad \text{for}\ z\in {_{\mathcal A} U}^+, f\in G_0^*(R) , [M]\in G_0(R).\end{equation}
$G_0^*(R)$ (resp. $G_0^*(R^\lambda)$) has the dual basis $\{\delta_M\mid [M]\in \mathbb B\}$ (resp. $\{\delta_{\mathcal M}\mid [\mathcal M]\in \mathbb B^\lambda\}$), where
$$\delta_{ M}([N])=\begin{cases} q^r \ \ \text{if}\ N\cong M\{r\}, \\ 0 \ \ \ \text{otherwise},\end{cases} \ \ \delta_{\mathcal M}([\mathcal N])=\begin{cases} q^r \ \ \text{if}\ \mathcal N\cong \mathcal M\{r\}, \\ 0 \ \ \ \text{otherwise}.\end{cases}$$
By \cite[Lemma 7.6(ii)]{LV2009}, $G_0^*(R)$ (resp. $G_0^*(R^\lambda)$) is generated by $\delta_{\mathbf 1}$ as ${_{\mathcal A} U}^+$-module, where $\mathbf 1=\K$ as a module over $R(0)=\K$ (resp. $R^\lambda(0)=\K$).

 \vskip 2mm
 Let $ {_{\mathcal A} V(\lambda)}^*$ be the graded dual of the $\mathcal A$-form $_{\mathcal A} V(\lambda)$. We define the ${_{\mathcal A} U}^+$-module structure of $ {_{\mathcal A} V(\lambda)}^*$ by the same way in (\ref{action}).  Then $ {_{\mathcal A} V(\lambda)}^*$  is generated by $\delta_{v_\lambda}$, the dual of the highest weight vector $v_\lambda$. Moreover, we have a ${_{\mathcal A} U}^+$-linear isomorphism
 $${_{\mathcal A} U}^+ \bigg/
(\sum_{i\in I^{\text{re}}}{_{\mathcal A} U}^+ e_i^{(\lambda_i+1)}+
\sum_{i\in I^\im \ \text{with}\ \lambda_i=0}{_{\mathcal A} U}^+e_i)\cong {_{\mathcal A} V(\lambda)}^*$$
which sends $\overline{1}$ to $\delta_{v_\lambda}$.
\vskip 2mm

\begin{lemma}\label{L7}
{\it We have
\begin{itemize}
\item[(i)] If $i\in I^\re$ and $m\geq \lambda_i+1$, then $e_i^{(m)}\cdot \delta_{\mathbf 1}=0$ in $G_0^*(R^\lambda)$.
\item[(ii)] If $i\in I^\im$ with $\lambda_i=0$, then $e_i\cdot \delta_{\mathbf 1}=0$ in $G_0^*(R^\lambda)$.
\end{itemize}
}
\begin{proof}
Under the assumption of (i),  we have $R^\lambda(mi)=0$ by the property of $V(i^m)$. For (ii), we have $R^\lambda(i)=0$ by the definition.
\end{proof}
\end{lemma}

 \vskip 2mm

By the above lemma, there are well-defined ${_{\mathcal A} U}^+$-linear epimorphisms
$${_{\mathcal A} U}^+\twoheadrightarrow G_0^*(R),\ \ {_{\mathcal A} V(\lambda)}^*\twoheadrightarrow G_0^*(R^\lambda),$$
which sends $1$ to $\delta_{\mathbf 1}$. Since $\mathbb B\cong B(\infty)$ and $\mathbb B^\lambda\cong  B(\lambda)$, we have for each $\nu$ that
$$\text{rank}({_{\mathcal A} U}^+_{\nu})=\text{rank}( G_0^*(R(\nu))),\ \ \text{rank}({_{\mathcal A} V(\lambda)}_{\lambda-\nu})=\text{rank}(  G_0^*(R^\lambda(\nu))).$$
Hence these maps must be isomorphisms. We obtain ${_{\mathcal A} U}^+$-linear isomorphisms
 $${_{\mathcal A} U}^+\cong G_0^*(R),\ \ {_{\mathcal A} V(\lambda)}^*\cong G_0^*(R^\lambda),$$
 and ${_{\mathcal A} V(\lambda)}\cong G_0(R^\lambda)$ by the duality.
 \vskip 6mm
\section{\textbf{Categorification of $V(\lambda)$}}
 \vskip 2mm
In this section, we construct an algebra $\mathscr R$ that is `lager' than $R$ in the sense that it has more irreducible classes when $i\in I^\iso$. We then form the cyclotomic quotient $\R^\lambda$ and realize $V(\lambda)$ $(\lambda\in P^+)$ as a subspace of $K_0(\R^\lambda)$, the Grothendieck group of $\R^\lambda$-\pMod, using the framework given in \cite{KK2012}.
\vskip 2mm
\subsection{The algebra  $\mathscr R$ }\
\vskip 2mm
Given a Borcherds-Cartan datum and a $\nu\in\N[I]$ with $\text{ht}{(\nu)}=n$. We define the generators of $\mathscr R(\nu)$ to be those of $R(\nu)$ and of the same degrees. They satisfy the  modified local relations:

\begin{align}
 \dcross{i}{j} \  =  \  \begin{cases} \quad \quad \quad \quad \quad \  0 & \text{ if } i= j\in I^\re, \\  \\ \   \Big( \dcrossL{-\frac{a_{ii}}{2}}{i}{i} \ + \ \dcrossR{-\frac{a_{ii}}{2}}{i}{i}\Big)^2  & \text{ if } i= j \ \text{and} \  (\alpha_i,\alpha_i) < 0,
                    \\ \\
                    \quad  \quad \quad \quad \ \ \dcrossA{i}{j} & \text{ if } \ (\alpha_i,\alpha_j)=0, \\
                    \\
                    \   \dcrossL{-a_{ij}}{i}{j} \ + \ \dcrossR{-a_{ji}}{i}{j}  & \text{ if } i\ne j \ \text{and} \  (\alpha_i,\alpha_j) \ne 0,
                  \end{cases}
\end{align}

\begin{equation}
\begin{aligned}
 \LU{{}}{i}{i} \   -  \ \RD{{}}{i}{i}\ =\ \dcrossA{i}{i} \quad\quad\quad \LD{{}}{i}{i}   \ - \ \RU{{}}{i}{i}\ =\ \dcrossA{i}{i}  \quad \text{ if } i \in I^\re,
\end{aligned}
\end{equation}

\begin{equation}
\begin{aligned}
 \LU{{}}{i}{j} \   =  \ \RD{{}}{i}{j}  \quad\quad\quad \LD{{}}{i}{j}   \ =  \ \RU{{}}{i}{j}  \quad \text{ otherwise},
\end{aligned}
\end{equation}

\begin{equation}
\quad \quad\quad\quad \quad\BraidL{i}{j}{i} \  -  \ \BraidR{i}{j}{i} \ =  {\sum_{c=0}^{-a_{ij}-1}} \ \threeDotStrands{i}{j}{i} \ \ \text{ if } i\in I^\re, i\ne j \ \text{and} \ (\alpha_i,\alpha_j) \ne 0,
 \end{equation}
 \begin{equation}
\BraidL{i}{j}{k} \  = \ \BraidR{i}{j}{k} \quad\quad \text{otherwise}.
\end{equation}

\vskip 2mm
 For $\ii\in \Seq(\nu)$, set $\mathscr{P}_{\ii}={\K}[x_1(\ii),\dots,x_n(\ii)]$ and form the $\K$-vector space $\mathscr{P}_{\nu}=\bigoplus_{\ii\in \Seq(\nu)}\mathscr{P}_{\ii}$. Each $\omega\in S_n$ acts on $\mathscr{P}_{\nu}$ by sending $x_a(\ii)$ to $x_{\omega (a)}(\omega(\ii))$.

 We assign a graph $\Lambda$ with vertices set $I$ and an edge between $i$ and $j$ if $i\neq j$ and  $(\alpha_i,\alpha_j) \ne 0$, and  choose an orientation for each edge. %Define the action of $\R({\nu})$ on $\Pol_{\nu}$ as follows: $_{\jj}\R(\nu)_{\ii}$ acts by $0$ on $\Pol_{\kk}$ if $\ii\neq\kk$.
 For $f\in \Pol_{\ii}$,  we set $1_{\ii}\cdot f=f$, $x_{k,\ii}\cdot f=x_k{(\ii)}f$ and
$$
\tau_{k,\ii}\cdot f=\begin{cases} \frac{f- s_kf}{x_k{(\ii)}-x_{k+1}(\ii)} & \text{if}\ i_k= i_{k+1}\in I^\re,\\
(x_k(\ii)^{-\frac{a_{ii}}{2}}+x_{k+1}(\ii)^{-\frac{a_{ii}}{2}})s_kf
& \text{if}\ i_k=i_{k+1}=i\in I^\im\backslash I^\iso,\\
{s_k}f & \text{if}\  (\alpha_{i_k},\alpha_{i_{k+1}})=0\ \text{or if}\ i_k\leftarrow i_{k+1}, \\
 (x_k(s_k\ii)^{-a_{i_{k+1}i_k}}+x_{k+1}(s_k\ii)^{-a_{i_ki_{k+1}}})s_kf & \text{if}\ i_k\rightarrow i_{k+1}.
\end{cases}
$$
It's easy to check that $\Pol_\nu$ is an $\R(\nu)$-module with the action above.
\vskip 2mm
Denote by $K_0(\R)$ ($\R=\bigoplus_{\nu}\R(\nu)$) the Grothendieck group of the category of finite generated gr-projective $\R$-modules. As in \cite{TW2023}, \cite{KL2009} and \cite{KL2011}, we  endow  $K_0(\R)$ with a twisted bialgebras structure and then we could obtain a twisted  bialgebras embedding $\Gamma: {_{\mathcal A} U}^-\hookrightarrow K_0(\R)$ given by
$$\begin{aligned}& f_i^{(n)}\mapsto P_{i^{(n)}}\quad \text{for}\ i\in I^\re,n\geq 0,\\ & f_i\mapsto P_{i}\quad \text{for}\ i\in I^\im.\end{aligned}$$
Here the meaning of $P_{i^{(n)}}$ and $P_i$ is the same as in section 1.2.
\vskip 2mm
\begin{remark}
$\Gamma$ is an isomorphism when $I^\iso=\emptyset$. If $I^\iso\neq \emptyset$, then $\Gamma$ is not surjective since for each $n$, the gr-irreducible classes of $\R(ni)$-modules can be labelled by the partitions of $n$.
\end{remark}

\vskip 2mm
\begin{lemma}\label{L8}
{\it Let $i\in I^\im$. Assume we add an additional relation $x_1^b=0$ in $\R(2i)$ for some $b\in \Z_{>0}$, then $x_2$ is nilpotent. }
\begin{proof}
If $i\in I^\iso$, then $x_2^b=\tau x_1^b\tau =0$. If $i\in I^\im\backslash I^\iso$, we set $a=-\frac{a_{ii}}{2}>0$. We show by downward induction on $0\leq k\leq 2b$ that $x_1^{ka}x_2^{b+2ab-ka}=0$. The case where $k=2b$ is obvious. If $ka\geq b$, the claim is also obvious. So we may assume that $ka< b$ and our claim is true for all $n> k$. Since $\tau^2=x_1^{2a}+2x_1^ax_2^a+x_2^{2a}$, we have
$$x_1^{ka}x_2^{b+2ab-ka}=x_1^{ka}x_2^{b+2ab-ka-2a}\tau^2-2x_1^{ka+a}x_2^{b+2ab-ka-a}-x_1^{ka+2a}x_2^{b+2ab-ka-2a}.$$
The last two terms in the right hand side are $0$ by the induction hypothesis.  It's enough to show that $2ab-ka-2a\geq 0$, but this follows from $2ab-ka-2a> 2a^2k-ka-2a\geq {-1}$. Then the lemma follows by taking $k=0$.
\end{proof}
\end{lemma}

\vskip 2mm

For each $i\in I$, define the functors
\begin{equation*}
\begin{aligned}
&E_{i}:\R(\nu+i)\text{-}\Mod\rightarrow  \R(\nu)\text{-}\Mod, \ N\mapsto 1_{\nu,i}  N, \\
&F_{i}:\R(\nu)\text{-}\Mod\rightarrow  \R(\nu+i)\text{-}\Mod, \ M\mapsto \Ind_{\nu,i}^{\nu+i} M\otimes \R(i),\\
&\overline{F}_{i}:\R(\nu)\text{-}\Mod\rightarrow  \R(\nu+i)\text{-}\Mod, \ M\mapsto \Ind_{i,\nu}^{\nu+i} \R(i)\otimes M,
\end{aligned}
\end{equation*}
The following proposition has been proved in \cite{KK2012}.

\vskip 2mm
\begin{proposition}\
{\it
\begin{itemize}
\item[(i)] There exists a natural isomorphism
\begin{equation*}
\begin{aligned}
E_iF_j\simeq\begin{cases}q_{i}^{-a_{ij}}F_jE_i& \text{if}\ i\neq j,\\
q_{i}^{-a_{ii}}F_iE_i\oplus (\rm{Id}\otimes\K[t_{i}]) & \text{if}\ i=j,
\end{cases}\\
\end{aligned}
\end{equation*}
where $t_{i}$ is an indeterminate of degree $2{r_{i}}$.
 %and $ \rm{Id}\otimes\boldsymbol{k}[t_{i}]:\rm{Mod}(R(\beta))\to \rm{Mod}(R(\beta))$ is the functor $M\to M\otimes\boldsymbol{k}[t_{i}]$.
 \item[(ii)] There exists a natural isomorphism $\overline{F}_j E_i\simeq E_i\overline{F}_j$ for $i\neq j$, and there is an exact sequence in $\R(\nu)\text{-}\Mod$
$$0\to \overline{F}_iE_iM\to E_i\overline{F}_iM\to q^{-(\nu,\alpha_{i})}M\otimes \K[t_{i}] \to 0.$$
\end{itemize}}
\end{proposition}

\vskip 2mm
Let $\lambda\in P^+$ and $i\in I$. As in Section 2, we define the cyclotomic algebra $\R^{\lambda}(\nu)$ to be the quotient of $\R(\nu)$ by the two sided ideal generated by $x_{1,\ii}^{\lambda_{i_1}}$ for all $\ii\in I$, and form
$$\R^\lambda=\bigoplus_{\nu}\R^{\lambda}(\nu), \ K_0(\R^\lambda)=\bigoplus_{\nu}K_0(\R^\lambda(\nu)),\ G_0(\R^\lambda)=\bigoplus_{\nu}G_0(\R^\lambda(\nu)).$$

\vskip 2mm

\begin{lemma}\label{L9}
{\it Let $\nu\in\N[I]$ with $\text{ht}(\nu)=n$. Then
 \begin{itemize}
\item [(i)] $x_{k,\ii}$ are nilpotent in $\R^\lambda(\nu)$ for all $\ii\in \Seq(\nu)$ and $1\leq k\leq n$. In particular, $\R^\lambda(\nu)$ is finite dimensional.
\item [(ii)] If $i \in I^\re$, then there exists $m\geq 0$ such that $\R^{\lambda}(\nu+ki)=0$ for any $k\geq m$.
\item [(iii)] If $i\in I^\im$ and $(\lambda-\nu)(h_i)=0$, then $\R^{\lambda}(\nu+i)=0$.
\end{itemize}}
\begin{proof}
The proof of (i) can be found in \cite[Proposition 2.3]{LV2009}. Here, Lemma \ref{L8} deals with some new cases  arise in the proof.

The proof of (ii) (resp. (iii)) is the same as \cite[Lemma 4.3(ii)]{KK2012} (resp. \cite[Lemma 4.4]{KKO2011}).
\end{proof}
\end{lemma}

\vskip 2mm

For each $i\in I$, define the functors
\begin{equation*}
\begin{aligned}
&E^\lambda_{i}:\R^\lambda(\nu+i)\text{-}\Mod\rightarrow  \R^\lambda(\nu)\text{-}\Mod, \ N\mapsto 1_{\nu,i}  N= 1_{\nu,i}\ \R^\lambda{(\nu+i)}\otimes_{\R^\lambda{(\nu+i)}}N, \\
&F^\lambda_{i}:\R^\lambda(\nu)\text{-}\Mod\rightarrow  \R^\lambda(\nu+i)\text{-}\Mod, \ M\mapsto \R^\lambda(\nu+i)1_{\nu,i}\otimes_{\R^\lambda(\nu)} M.
\end{aligned}
\end{equation*}
By a similar argument in \cite[Section 4]{KK2012}, one can prove:
\vskip 2mm
\begin{proposition}
{\it  The $\R^\lambda(\nu+i)1_{\nu,i}$ (resp. $1_{\nu,i}\R^\lambda(\nu+i)$) is a projective right (resp. left) $\R^\lambda(\nu)$-module. In particular, the functor $F^\lambda_{i}$ is exact, the functor $E^\lambda_{i}$  sends finitely generated projective modules to finitely generated projective modules.}
\end{proposition}

\vskip 2mm
\begin{remark}
Comparing with  \cite[Section 4]{KK2012}, we need to modify some of settings there to suit our case. First, for $\nu\in \N[I]$ with $\text{ht}(\nu)=n$, the intertwiner elements $g_k$ $(1\leq k\leq n)$ of $\R(\nu+i)$ defined in  \cite [(4.14)]{KK2012} should be changed into
$$g_{k}=\sum\limits_{\substack{{\ii \in \Seq(\nu+i)},\\ a_{i_k i_{k+1}}\leq 0}}\tau_{k,\ii}+\sum\limits_{\substack{\ii \in \Seq(\nu+i),\\ i_k=i_{k+1}\in I^\re}}\big(x_{k,\ii}-x_{k+1,\ii}-(x_{k,\ii}-x_{k+1,\ii})^2\tau_{k,\ii}\big).$$
Second, the elements $A,B\in \R(\nu+i)$ in \cite[Lemma 4.19]{KK2012} should be changed by
$$\begin{aligned} &A=\sum_{\ii\in \Seq (\nu)}\Big(x_{1,i\ii}^{\lambda_i}\prod_{\substack{1\leq k\leq n \\ i_k\neq i\ \text{and} \ a_{ii_k}\neq 0}}(x_{1,i\ii}^{-a_{ii_k}}+x_{k+1,i\ii}^{-a_{i_ki}})\prod_{\substack{1\leq k\leq n \\ i_k=i\in I^\im\backslash I^\iso}}(x_{1,i\ii}^{-\frac{a_{ii}}{2}}+x_{k+1,i\ii}^{-\frac{a_{ii}}{2}})^2\Big),\\
& B=\sum_{\ii\in \Seq (\nu)}\Big(x_{n+1,\ii i}^{\lambda_i}\prod_{\substack{1\leq k\leq n \\ i_k\neq i\ \text{and} \ a_{ii_k}\neq 0}}(x_{n+1,\ii i}^{-a_{ii_k}}+x_{k,\ii i}^{-a_{i_ki}})\prod_{\substack{1\leq k\leq n \\ i_k=i\in I^\im\backslash I^\iso}}(x_{n+1,\ii i}^{-\frac{a_{ii}}{2}}+x_{k,\ii i}^{-\frac{a_{ii}}{2}})^2\Big).\end{aligned}$$
\end{remark}

\vskip 2mm
By the exact same argument as in  \cite[Section 5]{KK2012}, one can prove:
\vskip 2mm
\begin{theorem}\
{\it
\begin{itemize}
\item[(1)] For $i\neq j\in I$, there exists a natural isomorphism
$$
 q_i^{-a_{ij}}F_j^{\lambda}E_i^{\lambda}\simeq E_i^{\lambda}F_j^{\lambda}.
$$
 \item[(2)] Let $\mu=\lambda-\nu$. Then there exists a natural isomorphisms of endofunctors on $\R^\lambda(\nu)\text{-}\Mod$:
\begin{itemize}
\item [(i)] If $\mu(h_i)\geq 0$, then
$$
\begin{aligned}
q_i^{-a_{ii}}F_i^{\lambda}E_i^{\lambda}\oplus\bigoplus_{k=0}^{\mu(h_i)-1}q_i^{2k}{\rm Id}\simeq E_i^{\lambda}F_i^{\lambda}.
\end{aligned}
$$
\item [(ii)]  If $\mu(h_i)< 0$, so $i\in I^\re$ and then
$$
\begin{aligned}
q_i^{-a_{ii}}F_i^{\lambda}E_i^{\lambda}\simeq E_i^{\lambda}F_i^{\lambda}\oplus\bigoplus_{k=0}^{-\mu(h_i)-1}q_i^{-2k-2}{\rm Id}.
\end{aligned}
$$
\end{itemize}\end{itemize}}
\end{theorem}

\vskip 2mm
Define the functors $\mathcal E_i^\lambda, \mathcal F_i^\lambda, K_i$ on each $\R^\lambda(\nu)\text{-}\Mod$ by
$$\mathcal E_i^\lambda=E_i^\lambda,\ \mathcal F_i^\lambda=q_i^{1-(\lambda-\nu)(h_i)}F_i^\lambda, \ K_i=q_i^{(\lambda-\nu)(h_i)}$$
Then by the above theorem, the induced endomorphisms on $K_0(\R^\lambda)$ and $G_0(\R^\lambda)$ satisfy
$$[\mathcal E_i^\lambda, \mathcal F_i^\lambda]=\delta_{ij}\frac{K_i-K_i^{-1}}{q_i- q_i^{-1}} $$
for all $i,j \in I$. Let $K_0(\R^\lambda)_{\Q(q)}=\Q(q)\otimes_{\Z[q,q^{-1}]}K_0(\R^\lambda)$ and $G_0(\R^\lambda)_{\Q(q)}=\Q(q)\otimes_{\Z[q,q^{-1}]}G_0(\R^\lambda)$. Then Lemma \ref{L9} and \cite[Proposition B.1]{KMPY1996} yield $K_0(\R^\lambda)_{\Q(q)}$ and $G_0(\R^\lambda)_{\Q(q)}$ integrable  $U_q(\g)$-modules.

Taking the subspace $V^\lambda$ of $K_0(\R^\lambda)_{\Q(q)}$ which is spanned by the elements of the form $ F_{i_1}^\lambda\cdots F_{i_s}^\lambda\cdot\mathbf 1_\lambda$ for $s\geq 0, i_1\dots i_s\in I$, where $\mathbf 1_\lambda$ is the trivial module over $\R^\lambda(0)$. Then  $V^\lambda$ is an integrable highest weight module with highest weight $\lambda$. So we have
 $V^\lambda \cong V(\lambda)$. The subspace of $G_0(\R^\lambda)_{\Q(q)}$ spanned by the same elements is also isomorphic to
 $V(\lambda)$.

\vskip 2mm
\begin{remark}
If $I^\iso= \emptyset$, then $V^\lambda=K_0(\R^\lambda)_{\Q(q)}$ and the Grothendieck group $K_0(\R^\lambda)$ is isomorphic to the $\mathcal A$-form ${_{\mathcal A} V(\lambda)}$. This follows from a similar argument in the last part of \cite[Section 6]{KK2012}.

Although our construction of $\R$  is  more concise, we do not get a stronger result than \cite{KKO2011}.
\end{remark}

\vskip 6mm

\bibliographystyle{amsplain}

%%%%%%%%
\end{document}